\documentclass[a4paper,12pt]{article}
\usepackage{amsmath}
\usepackage{amsfonts}
\usepackage{amssymb}
\usepackage{latexsym}
\usepackage{epsfig}
\usepackage{graphicx}
\usepackage{oldgerm}
\usepackage{theorem}

\def\R{\mathbb{R}}
\def\N{\mathbb{N}}
\def\H{\mathbb{H}}
\def\C{\mathbb{C}}
\def\W{\mathbb{W}}
\def\Z{\mathbb{Z}}

\def\I{\mathbb{I}}

\def\rP{{\rm P}}
\def\rE{{\rm E}}
\def\P{\mathbb{P}}
\def\E{\mathbb{E}}
\def\bE{{\bf E}}

\def\1{{\bf 1}}
\def\0{{\bf 0}}

\def\F{{\cal F}}
\def\S{{\cal S}}

\def\cH{{\cal H}}

\def\x{\mib{x}}
\def\y{\mib{y}}
\def\z{\mib{z}}
\def\X{\mib{X}}
\def\Y{\mib{Y}}
\def\B{\mib{B}}
\def\bZ{\mib{Z}}
\def\t{\mib{t}}
\def\f{\mib{f}}
\def\a{\mib{a}}
\def\b{\mib{b}}
\def\vlambda{\mib{\lambda}}
\def\valpha{\mib{\alpha}}

\def\d={\stackrel{\rm d}{=}}

\def\mM{\mathfrak{M}}
\def\mbK{\mathbb{K}}

\def\supp{{\rm supp}\ }
\def\mbK{\mathbb{K}}

\def\c{{\rm C}}
\def\Det{\mathop{\mathrm{Det}}}

\def\LE{{\rm LE}}
\theorembodyfont{\itshape}
\newtheorem{thm}{Theorem}[section]

\newcommand{\mib}[1]{\mbox{\boldmath $#1$}}
\newcommand{\SSC}[1]{\section{#1}\setcounter{equation}{0}}

\begin{document}

\title{\bf
Bessel process, Schramm-Loewner evolution,
and Dyson model \\
\vskip 0.3cm
{\large -- Complex Analysis applied to\\
Stochastic Processes and Statistical Mechanics --}
\footnote{
This manuscript is prepared for the proceedings of
the 9th Oka symposium, held at Nara Women's University,
4-5 December 2010.
}}
\author{
Makoto Katori
\footnote{
Department of Physics,
Faculty of Science and Engineering,
Chuo University, 
Kasuga, Bunkyo-ku, Tokyo 112-8551, Japan;
e-mail: katori@phys.chuo-u.ac.jp
}}
\date{24 March 2011}
\pagestyle{plain}
\maketitle
\begin{abstract}
{\it Bessel process} is defined as the radial part of 
the Brownian motion (BM) in the $D$-dimensional space,
and is considered as a one-parameter family
of one-dimensional diffusion processes 
indexed by $D$, BES$^{(D)}$.
First we give a brief review of BES$^{(D)}$,
in which $D$ is extended to be a continuous
positive parameter.
It is well-known that $D_{\rm c}=2$ is 
the critical dimension such that,
when $D \geq D_{\rm c}$ (resp. $D < D_{\rm c}$),
the process is transient (resp. recurrent).
{\it Bessel flow} is a notion such that 
we regard BES$^{(D)}$ with a fixed $D$ as 
a one-parameter family of initial value $x > 0$.
There is another critical dimension
$\overline{D}_{\rm c}=3/2$ and, in the intermediate
values of $D$,
$\overline{D}_{\rm c} < D < D_{\rm c}$,
behavior of Bessel flow is highly nontrivial.
The dimension $D=3$ is special, since
in addition to the aspect that BES$^{(3)}$ is 
a radial part of the three-dimensional BM,
it has another aspect as
a {\it conditional BM to stay positive}.

Two topics in probability theory
and statistical mechanics,
the {\it Schramm-Loewner evolution} (SLE)
and the {\it Dyson model} 
({\it i.e.}, Dyson's BM model with parameter $\beta=2$),
are discussed.
The SLE$^{(D)}$ is introduced as a 
`complexification' of Bessel flow on the
upper-half complex-plane,
which is indexed by $D > 1$.
It is explained that the existence of
two critical dimensions $D_{\rm c}$ and $\overline{D}_{\rm c}$
for BES$^{(D)}$ makes SLE$^{(D)}$
have three phases;
when $D \geq D_{\rm c}$ the SLE$^{(D)}$ path is simple,
when $\overline{D}_{\rm c} < D < D_{\rm c}$
it is self-intersecting but not dense,
and when $1 < D \leq \overline{D}_{\rm c}$
it is space-filling.
The Dyson model is introduced as a multivariate extension
of BES$^{(3)}$. 
By `inheritance' from BES$^{(3)}$,
the Dyson model has two aspects;
(i) as an eigenvalue process of a Hermitian-matrix-valued
BM, and (ii) as a system of BMs conditioned never to
collide with each other,
which we simply call the {\it noncolliding BM}.
The noncolliding BM is constructed as a harmonic
transform of absorbing BM in the Weyl chamber of type A,
and as a complexification of this construction,
the {\it complex BM representation} is proposed
for the Dyson model.
Determinantal expressions for spatio-temporal 
correlation functions with the {\it asymmetric correlation kernel
of Eynard-Mehta type} are direct consequence of this
representation.
In summary, `parenthood' of BES$^{(D)}$ and SLE$^{(D)}$,
and that of BES$^{(3)}$ and the Dyson model are clarified.

Other related topics concerning extreme value distributions
of noncolliding diffusion processes,
statistics of characteristic polynomials of random matrices,
and scaling limit of Fomin's determinant for 
loop-erased random walks are also given.

We note that the name of Bessel process is
due to the {\it special function} called
the modified Bessel function,
SLE is a stochastic time-evolution of
{\it complex analytic function}
(conformal transformation),
and the Weierstrass canonical product representation of
{\it entire functions} plays an important role
for the Dyson model.
Complex analysis is effectively applied to
study stochastic processes of interacting particles and
statistical mechanics models exhibiting
critical phenomena and fractal structures
in equilibrium and nonequilibrium states. \\
{\bf Keywords} \,
Complexification, Multivariate extension,
Conformal transformation, 
Random matrices, 
Entire functions
\end{abstract}

\tableofcontents
\vspace{3mm}

\SSC{Family of Bessel processes}
\subsection{One-dimensional and $D$-dimensional
Brownian motions}

We consider motion of a Brownian particle in one dimensional 
space $\R$ starting from $x \in \R$ at time $t=0$.
At each time $t>0$ particle position is randomly distributed, 
and each realization of path is labeled by a parameter $\omega$.
Let $\Omega$ be the sample path space and
$B^{x}(t, \omega)$ denote the position of the 
Brownian particle at time $t>0$, whose path is
realized as $\omega \in \Omega$.
Let $(\Omega, \F, \rP)$ be a probability space.
The {\it one-dimensional standard Brownian motion} (BM),
$\{B^{x}(t, \omega)\}_{t \in [0, \infty)}, x \in \R$,
has the following three properties.

\begin{description}
\item{1.} \quad
$B^{x}(0,w)=x$ with probability one
(abbr. w.p.1).

\item{2.} \quad
For any fixed $\omega \in \Omega$,
$B^{x}(t, \omega)$ is a real continuous function of $t$.
In other words, $B^{x}(t)$ has a continuous path.

\item{3.} \quad
For any sequence of times, 
$t_0 \equiv 0 < t_1< \cdots < t_M, 
M \in \N \equiv \{1,2,3, \dots\}$,
the increments $\{B^{x}(t_i)-B^{x}(t_{i-1})\}_{i=1}^{M}$
are independent, 
and distribution of each increment is normal with mean
$m=0$ and variance $\sigma^2=t_{i}-t_{i-1}$.
It means that for any $1 \leq i \leq M$ and $\ell < r$,
$$
\rP(B^{x}(t_i)-B^{x}(t_{i-1}) \in [\ell, r])
=\int_{\ell}^{r} p_{t_i-t_{i-1}}(\delta|0) d \delta,
$$
where we define for $a,b \in \R$
\begin{equation}
p_t(b|a)= \left\{
\begin{array}{ll}
\displaystyle{\frac{1}{\sqrt{2 \pi t}} e^{-(a-b)^2/2t}},
& \quad \mbox{for $t > 0$},
\cr
\delta(a-b),
& \quad \mbox{for $t=0$}.
\end{array} \right.
\label{eqn:pt1}
\end{equation}
\end{description}

If we write the conditional probability as
$\rP(\cdot | {\cal C})$,
where ${\cal C}$ denotes the condition,
the third property given above implies that
for any $0 \leq s \leq t$
$$
\rP(B^{x}(t) \in A | B^{x}(s)=a)
=\int_{A} p_{t-s}(b|a) db
$$
holds $^{\forall} A \subset \R, ^{\forall} a \in \R$.
Then the probability that the BM is observed
in a region $A_i \subset \R$ at time $t_i$
for each $i=1, 2, \dots, M$ is given by
\begin{equation}
\rP(B^{x}(t_i) \in A_i, i=1,2, \dots, M)
=\int_{A_1}dx_1 \cdots \int_{A_M} dx_M \,
\prod_{i=1}^{M} p_{t_i-t_{i-1}}(x_i|x_{i-1}),
\label{eqn:BM3}
\end{equation}
where $x_0 \equiv x$.
The formula (\ref{eqn:BM3}) means that for any fixed $s \geq 0$,
under the condition that $B^{x}(s)$ is given,
$\{B^{x}(t): t \leq s\}$ and $\{B^{x}(t) : t > s\}$ 
are independent.
This independence of the events in the future
and those in the past is called {\it Markov property}
\footnote{
A positive random variable $\tau$ is called {\it Markov time}
if the event $\{\tau \le u\}$ is determined 
by the behavior of the process
until time $u$ and independent of that after $u$.
The Brownian motion satisfies the property obtained by changing
any deterministic time $s>0$ into any Markov time $\tau$ 
in the definition of Markov property given here.
It is called a {\it strong Markov property}.
A stochastic process which is strong Markov and has a continuous path
almost surely is called a {\it diffusion process}.
}.
The integral kernel $p_{t}(y|x)$ is called
the {\it transition probability density function}
of the BM.
As defined by (\ref{eqn:pt1}),
$p_{t}(y|x)$ is nothing but the probability density function
of the normal distribution (the Gaussian distribution)
with mean $m=x$ and variance $\sigma^2=t$.
It should be noted that $p_{t}(y|x)=p_{t}(x|y)$
and $u_{t}(x) \equiv p_{t}(y|x)$
is a unique solution of the {\it heat equation} 
({\it diffusion equation})
\begin{equation}
\frac{\partial}{\partial t} u_t(x)
=\frac{1}{2} \frac{\partial^2}{\partial x^2} u_t(x),
\quad x \in \R, \quad t \in [0, \infty)
\label{eqn:heat}
\end{equation}
with the initial condition $u_0(x)=\delta(x-y)$.
Therefore, $p_{t}(y|x)$ is also called
the {\it heat kernel}.

Let $D \in \N$ denote the spatial dimension.
For $D \geq 2$, the $D$-dimensional BM in $\R^D$ starting
from the position $\x=(x_1, \dots, x_D) \in \R^D$
is defined by the following $D$-dimensional
vector-valued BM,
\begin{equation}
\B^{\x}(t)=(B^{x_1}_1(t), B^{x_2}_2(t), \dots,
B^{x_D}_D(t)),
\quad t \in [0, \infty),
\label{eqn:DBM}
\end{equation}
where $\{B^{x_i}_i(t)\}_{i=1}^{D}$
are independent one-dimensional standard BMs.

\subsection{$D$-dimensional Bessel process}

\begin{figure}
\begin{center}
\includegraphics[width=0.5\linewidth]{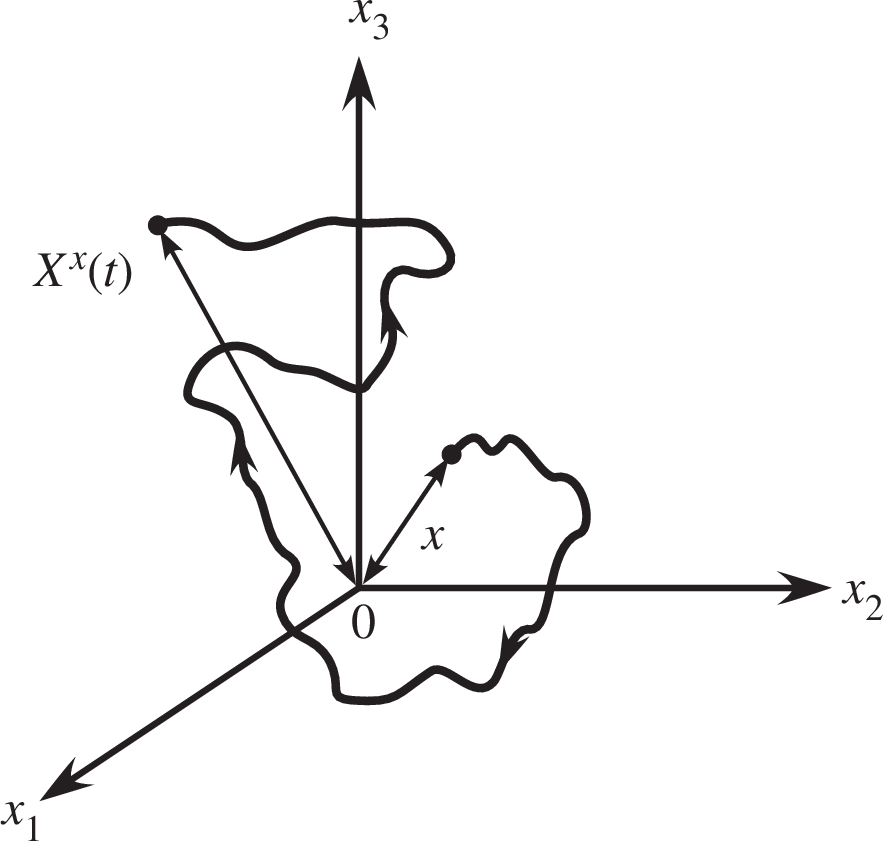}
\end{center}
\caption{\small
The $D$-dimensional Bessel process $X^{x}(t)$
is defined as the radial part of the BM in the
$D$-dimensional space.
The initial value $x$ of the Bessel process
is the distance between the origin and 
the position from which the BM is started.}
\label{fig:BESD}
\end{figure}

For $D \in \N$, the $D$-dimensional {\it Bessel process}
is defined as the absolute value
({\it i.e.}, the radial coordinate) of the
$D$-dimensional BM,
\begin{eqnarray}
X^{x}(t) &\equiv& |\B^{\x}(t)|
\nonumber\\
&=& \sqrt{ B^{x_1}_1(t)^2+ \cdots
+ B^{x_D}_D(t)^2},
\quad t \in [0, \infty),
\label{eqn:BES1}
\end{eqnarray}
where the initial value is given by
$X^{x}(0)=x=|\x|=\sqrt{x_1^2+ \dots+ x_D^2} \geq 0$.
By definition $X^{x}(t)$ is nonnegative,
$X^{x}(t) \in \R_{+} \equiv \{x \in \R: x \geq 0\}$.
See Fig. \ref{fig:BESD}.

By this definition, $X^{x}(t)$ is a functional of
$D$-tuples of stochastic processes 
$\{B^{x_i}_i(t)\}_{i=1}^{D}$.
In order to describe the statistics of a function 
of several random variables, we have to see 
`propagation of error'.
For stochastic processes, by {\it It\^o's formula}
we can readily obtain an equation for the stochastic 
process that is defined as a functional of several stochastic
processes.
In the present case, we have the following equation,
\begin{equation}
dX^x(t)=dB(t)+ \frac{D-1}{2} \frac{dt}{X^x(t)},
\quad t \in [0, \infty), \quad x>0.
\label{eqn:BES2}
\end{equation}
The first term of the RHS, $dB(t)$, denotes the infinitesimal
increment of a one-dimensional standard BM starting
from the origin at time $t=0$, $B(t)=B^0(t)$.
It should be noted that $B(t)$ is a different BM
from any $B^{x_i}_i(t), 1 \leq i \leq D$,
which were used to define $X^{x}(t)$ in Eq.(\ref{eqn:BES1}).
Here we assume $X^{x}(t) > 0$.
Then, if $D >1$, for an infinitesimal increment of
time $dt > 0$, the second term in the RHS of (\ref{eqn:BES2})
is positive.
It means that there is a drift to increase the value of $X^{x}(t)$.
This drift term is increasing in $D$ and decreasing in $X^{x}(t)$.
Since as $X^{x}(t) \searrow 0$, the drift term
$\nearrow \infty$, it seems that 
a `repulsive force' is acting to the $D$-dimensional BM, 
$\B^{\x}(t), |\x| > 0$ to keep the distance
from the origin positive, $X^{x}(t)=|\B^{\x}(t)| > 0$,
in order to avoid collision of the Brownian particle
at the origin.
Such a differential equation as (\ref{eqn:BES2}), which involves
random fluctuation term and drift term is called
a {\it stochastic differential equation} (SDE).

What is the origin of the repulsive force between 
the $D$-dimensional BM and the origin?
Why $\B^{\x}(t)$ starting from a point $\x \not=0$
does not want to return to the origin ?
Why the strength of the outward drift is increasing
in the dimension $D > 1$ ?

There is no positive reason for $\B^{\x}(t)$
to avoid visiting the origin, since by definition 
(\ref{eqn:DBM}) all components $B^{x_i}_i(t)$ enjoy
independent BMs.
As the dimension of space $D$ increases, however,
the possibility {\it not} to visit the origin
(or the fixed special point) increases, since among $D$ directions
in the space only one direction is toward the origin
(or the fixed special point) and other $D-1$ directions
are orthogonal to it.
If one know the second law of thermodynamics,
which is also called the {\it law of increasing entropy},
one will understand that we would like to say here that
the repulsive force acting from the origin to the Bessel process
is an `entropy force'.
(Note that the physical dimension of entropy 
[J/K] is different from that of force [J/m].)
Anyway, the important fact is that, while the variance
(quadratic variation)
of the standard BM is fixed as $(dB(t))^2=dt$
for a given $dt >0$,
the strength of repulsive drift is 
increasing in $D$.
Then, the return probability of $X^{x}(t), x >0$
to the origin should be a decreasing function of $D$.

Let $p^{(D)}_t(y|x)$ be the transition probability density of
the $D$-dimensional Bessel process.
We can show that, for any $y \in \R_{+}$,
$u^{(D)}_t(x) \equiv p^{(D)}_t(y|x), x>0$ solves the
following partial differential equation (PDE)
\begin{equation}
\frac{\partial}{\partial t} u^{(D)}_t(x)
=\frac{1}{2} \frac{\partial^2}{\partial x^2} u^{(D)}_t(x)
+\frac{D-1}{2x}
\frac{\partial}{\partial x} u^{(D)}_t(x)
\label{eqn:BES3}
\end{equation}
under the initial condition $u_0^{(D)}(x)=\delta(x-y)$,
which is called the {\it backward Kolmogorov equation}
for the $D$-dimensional Bessel process.
We can see clear correspondence between the SDE (\ref{eqn:BES2})
and the PDE (\ref{eqn:BES3}).
As shown by (\ref{eqn:heat}),
the BM term, $dB(t)$, in (\ref{eqn:BES2})
is mapped to the diffusion term 
$(1/2) \partial^2 u^{(D)}_t(x)/\partial x^2$
in (\ref{eqn:BES3}).
In (\ref{eqn:BES3}) the drift term is given by using
the spatial derivative $\partial/\partial x$
representing the outward drift with the coefficient
$(D-1)/2x$ corresponding to the factor
$(D-1)/(2X^{x}(t))$ of the second term in (\ref{eqn:BES2}).
The solution is given by
\begin{equation}
p^{(D)}_t(y|x)
=  \left\{ \begin{array}{ll}
\displaystyle{
\frac{1}{t} \frac{y^{\nu+1}}{x^{\nu}}
e^{-(x^2+y^2)/2t}
I_{\nu} \left( \frac{xy}{t} \right)},
& \quad t>0, x >0, y \geq 0, \cr
\displaystyle{
\frac{y^{2\nu+1}}{2^{\nu} t^{\nu+1} \Gamma(\nu+1)} e^{-y^2/2t}},
& \quad t>0, x=0, y \geq 0 \cr
& \cr
\delta(y-x),
& \quad t=0, x, y \geq 0,
\end{array} \right.
\label{eqn:pD}
\end{equation}
where $I_{\nu}(z)$ is the {\it modified Bessel function}
of the first kind defined by 
\begin{equation}
I_{\nu}(z) = \sum_{n=0}^{\infty} 
\frac{1}{\Gamma(n+1) \Gamma(n+1+\nu)}
\left( \frac{z}{2} \right)^{2n+\nu}
\label{eqn:I}
\end{equation}
with the gamma function 
$\Gamma(z)=\int_0^{\infty} e^{-u} u^{z-1} du$,
and the index $\nu$ is specified by the dimension $D$ as
\begin{equation}
\nu=\frac{D-2}{2}
\quad \Longleftrightarrow \quad
D=2(\nu+1).
\label{eqn:nuD}
\end{equation}
This fact that $p^{(D)}_t(y|x)$ is expressed by
using $I_{\nu}(z)$ gives the reason why
the process $X^{x}(t)$ is called the Bessel process.

\begin{figure}
\begin{center}
\includegraphics[width=1.0\linewidth]{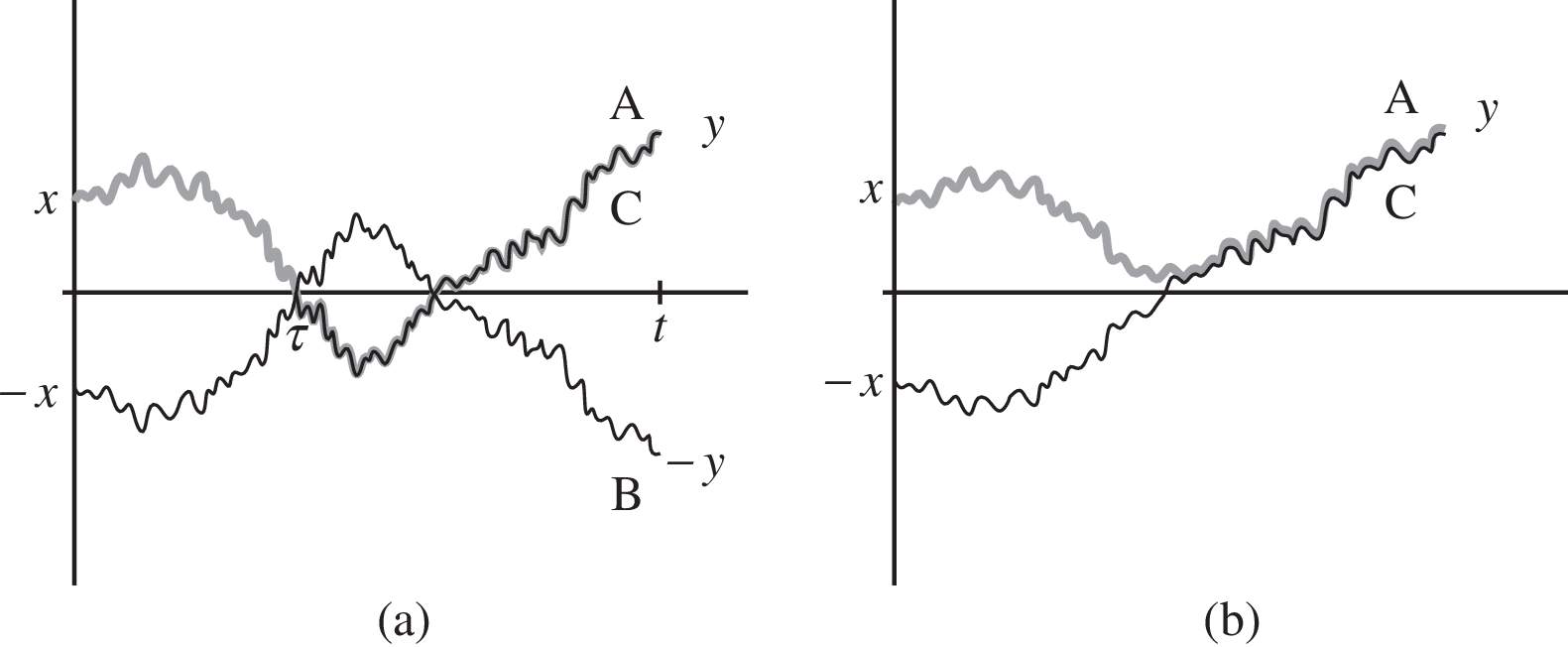}
\end{center}
\caption{\small
(a) One realization of Brownian path
from $x>0$ to $y>0$ is drawn (path A),
which visits the nonpositive region
$\R_{-} \equiv \{x \in \R: x \leq 0 \}$.
The path B is a mirror image of path A with respect
to the origin $x=0$, which is running from $-x<0$
to $-y < 0$. 
The first time when the path B hits the origin is
denoted by $\tau$.
The path C is a combination of a part of path B
up to time $\tau$ 
and a part of path A after $\tau$
such that it runs from $-x < 0$ to $y>0$.
There establishes bijection between path A and path C,
which have the same weight as Brownian paths.
Since the Brownian path A contributes
to $p_t(y|x)$ and the Brownian path B does to $p_t(y|-x)$,
such a path from $x>0$ to $y>0$ visiting $\R_-$ is cancelled
in $p_{t}^{\rm abs}(y|x)$.
(b) Each path, which does not visit $\R_-$,
gives positive contribution to $p_t^{\rm abs}(y|x)$,
since in this case the weight of Brownian path A
contributing to $p_{t}(y|x)$ is bigger than
that of Brownian path C contributing $p_{t}(y|-x)$.
In summary, $p_{t}^{\rm abs}(y|x)=p_t(y|x)-p_t(y|-x)$
gives the total weight of Brownian paths, 
which do not hit the origin and thus are not absorbed
at the wall at the origin.}
\label{fig:absBM}
\end{figure}
\begin{figure}
\begin{center}
\includegraphics[width=0.6\linewidth]{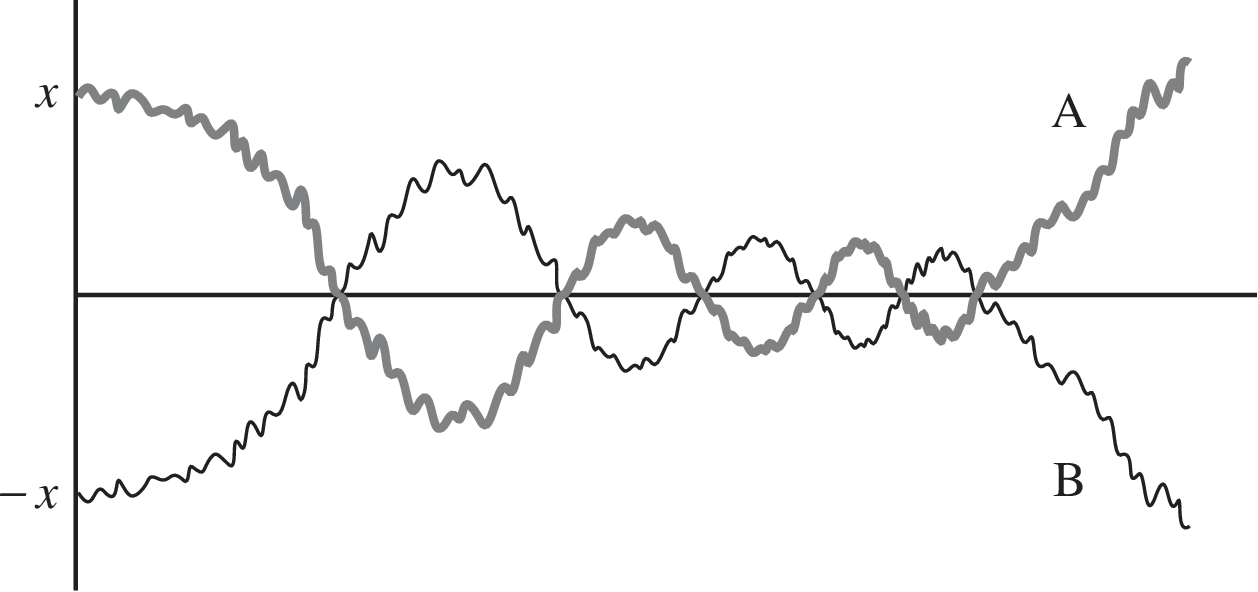}
\end{center}
\caption{\small
The one-dimensional BM visits the origin frequently.
For a Brownian path starting from $x>0$ (path A),
its mirror image with respect to the origin is drawn (path B),
which starts from $-x<0$.
If we observe the motion only in the nonnegative region $\R_+$,
the superposition of Brownian paths A and B gives
a path of a reflecting BM, where a reflecting wall
is put at the origin.}
\label{fig:refBM}
\end{figure}
When $D=3$, $\nu=1/2$ by (\ref{eqn:nuD}),
and we can use the equality 
$I_{1/2}(z)=\sqrt{2/\pi z} \sinh z
=(e^{z}-e^{-z})/\sqrt{2 \pi z}$.
Then (\ref{eqn:pD}) gives
\begin{equation}
p_t^{(3)}(y|x)=\frac{y}{x} \Big\{
p_t(y|x)-p_t(y|-x) \Big\}
\label{eqn:p3}
\end{equation}
for $t >0, x >0, y \geq 0$, where $p_t(y|x)$ is 
the transition probability density (\ref{eqn:pt1}) of BM.
If we put 
$p^{\rm abs}_t(y|x)=p_t(y|x)-p_t(y|-x)$, we see
$p^{\rm abs}_t(0|x)=0$ for any $x > 0$,
since BM is a symmetric process.

As shown by Fig.\ref{fig:absBM}, $p^{\rm abs}_t(y|x), x, y \geq 0$,
gives the transition probability density of the absorbing BM,
in which an absorbing wall is put at the origin and,
if the Brownian particle starting from $x>0$ arrives at the origin,
it is absorbed there and the motion is stopped.
By absorption, the total mass of paths from $x >0$ to $y>0$
is then reduced.
The factor $y/x$ appearing in (\ref{eqn:p3}) is for 
renormalization so that
$\int_{\R_+} p^{(3)}_t(y|x)dy=1, \forall t >0, \forall x >0$.
We regard this renormalization procedure from $p^{\rm abs}_t$
to $p^{(3)}_t$ as a transformation.
Since $x$ is a one-dimensional harmonic function 
in a rather trivial sense
$\Delta^{(1)} x \equiv d^2 x/d x^2=0$,
we say that the three-dimensional Bessel process is
an {\it harmonic transform}
($h$-transform) of the one-dimensional absorbing BM
in the sense of Doob \cite{Doo84}.
This implies the equivalence between the three-dimensional
Bessel process and `the one-dimensional BM
{\it conditioned to stay positive}'.
We will discuss such equivalence of processes in Section 2
more detail.
Here we put emphasize the fact that $p^{(3)}_t(0|x)=0,
\forall x >0$.
It means that the three-dimensional Bessel process does not
visit the origin.
When $D=3$, the outward drift is strong enough to avoid
any visit to the origin.
Moreover, we can prove that for any $x>0$, $X^{x}(t) \to \infty$
as $t \to \infty$ w.p.1
and we say the process is {\it transient}.

When $D=1$, $\nu=-1/2$ by (\ref{eqn:nuD}) and we use
the equality
$I_{-1/2}(z)=\sqrt{2/\pi z} \cosh z
=(e^z+e^{-z})/\sqrt{2 \pi z}$.
In this case (\ref{eqn:pD}) gives
\begin{equation}
p^{(1)}_t(y|x)=p_t(y|x)+p_t(y|-x)
\label{eqn:p1}
\end{equation}
for $t >0, x, y \geq 0$.
As shown by Fig.\ref{fig:refBM}, (\ref{eqn:p1}) means the
equivalence between the one-dimensional Bessel process
and `the one-dimensional BM with a reflecting wall at the origin'.
This is of course a direct consequence of the
definition of Bessel process (\ref{eqn:BES1}),
since it gives $X^{x}(t)=|B^{x}(t)|$ in $D=1$.
The important fact is that the one-dimensional BM starting
from $x \not=0$ visits the origin frequently
and we say that the one-dimensional Bessel process is 
{\it recurrent}.
(Remark that in Eqs. (\ref{eqn:BES2}) and (\ref{eqn:BES3}),
the drift terms vanish when $D=1$.
So we have to assume the reflecting boundary condition
at the origin when we discuss the one-dimensional Bessel process
instead of the one-dimensional BM.)

Now the following question is addressed:
At which dimension the Bessel process changes its property
from recurrent to transient ?

Before answering this question, here we would like to extend
the setting of the question.
Originally, the Bessel process was defined by (\ref{eqn:BES1})
for $D \in \N$.
We find that, however, the modified Bessel function 
(\ref{eqn:I}) is an analytic function of 
$\nu$ for all values of $\nu$.
So we will be able to define the Bessel process
for any positive value of dimension $D>0$
as the diffusion process in $\R_+$
such that the transition probability density function
is given by (\ref{eqn:pD}), where the index 
$\nu > -1$ is determined by (\ref{eqn:nuD})
for each value of $D$.
(In the SDE, (\ref{eqn:BES2}), we assume the
reflecting boundary condition at the origin
for $0 < D < 2$.)
Now we introduce an abbreviation BES$^{(D)}$
for the $D$-dimensional Bessel process,
$D > 0$
\footnote{
Another characterization of BES$^{(D)}$ for
fractional dimensions $D$ is given by the following.
Let $\nu \in \R$ and consider a BM with a constant
drift $\nu$, 
$B^{y}(t)+\nu t$, which starts from $y \in \R$ at time $t=0$.
The {\it geometric BM} with drift $\nu$ is defined as
$\exp(B^{y}(t)+\nu t), t \geq 0$.
For each $t \geq 0$, if we define random time change 
$ t \mapsto A_t$ by
$$
A_t=\int_0^{t} \exp \{ 2(B_s+\nu s)\} ds,
$$
then the following relation is established,
$$
X^{x}(A_t)=\exp(B_t^x+\nu t), \quad t \geq 0,
$$
where $X^{x}(A_t)$ is the BES$^{(D)}$ with
$D=2(\nu+1)$ at time $A_t$ starting from $x=e^y$.
The above formula is called {\it Lamperti's relation}
\cite{Lam72,Yor01}
}.

For BES$^{(D)}$ starting from $x>0$, denote
its first visiting time at the origin by
\begin{equation}
T^x=\inf \{t > 0: X^x(t)=0 \}.
\label{eqn:Tx}
\end{equation}
The answer of the above question is given by
the following theorem.

\begin{thm}
\label{thm:BESD1}
\begin{description}
\item{\rm (i)} \quad
$D \geq 2 \quad \Longrightarrow \quad
T^{x}=\infty, ^{\forall}x >0$, w.p.1.
\item{\rm (ii)} \quad
$D > 2 \quad \Longrightarrow \quad
\displaystyle{\lim_{t \to \infty}} X^{x}(t)=\infty,
\, ^{\forall}x > 0$, w.p.1,
{\it i.e.} the process is transient.
\item{\rm (iii)} \quad
$D=2 \quad \Longrightarrow \quad
\displaystyle{\inf_{t >0}} \, X^{x}(t)=0, \, 
^{\forall}x >0$, w.p.1D \\
That is, BES$^{(2)}$ starting from $x>0$ does not
visit the origin, but it can visit any neighbor of
the origin.
\item{\rm (iv)} \quad
$D < 2 \quad \Longrightarrow \quad
T^x < \infty, \, ^{\forall} x > 0$, w.p.1,
{\it i.e.} the process is recurrent.
\end{description}
\end{thm}
\vskip 0.5cm

\subsection{Bessel flow and Cardy's formula}

In the previous subsection, we have defined
the BES$^{(D)}$ for positive continuous values of 
dimension $D >0$ and studied dependence of the probability law
of process on $D$.
Theorem \ref{thm:BESD1} states that the two-dimension is 
a {\it critical dimension},
$$
 D_{\rm c}=2,
$$
for competition between the two effects acting 
the Bessel process,
the `random force' (the martingale term)
and the `entropy force' (the drift term)
in (\ref{eqn:BES2}) and (\ref{eqn:BES3}):
when $D > D_{\rm c}$, the latter dominates the former
and the process becomes transient, and when 
$D < D_{\rm c}$, the former is relevant and recurrence
to the origin of the process is realized frequently.

Here we show that there is another critical dimension,
$$
\overline{D}_{\rm c}=\frac{3}{2}.
$$
In order to characterize the transition at $\overline{D}_{\rm c}$,
we have to investigate dependence of the behavior
of $X^{x}(t)$ on initial value $x >0$.
We call the one-parameter family 
$\{X^{x}(t)\}_{x >0}$ the {\it Bessel flow} 
for each fixed $D>0$.

For $0 < x < y$, we trace the motions of two BES$^{(D)}$'s
starting from $x$ and $y$ by solving
(\ref{eqn:BES2}) using the {\it common} BM, $B(t), t \geq 0$,
\begin{eqnarray}
&& X^{x}(t)= x+ B(t)+\frac{D-1}{2}
\int_0^{t} \frac{ds}{X^{x}(s)},
\nonumber\\
&& X^{y}(t)= y+ B(t)+\frac{D-1}{2}
\int_0^{t} \frac{ds}{X^{x}(s)}, \quad t \geq 0.
\nonumber
\end{eqnarray}
By considering the coupling of the two processes,
we can show that
\begin{eqnarray}
 x<y \quad &\Longrightarrow& \quad
X^x(t) < X^{y}(t), \quad t < T^{x}
\quad \mbox{w.p.1}
\nonumber\\
&\Longrightarrow& \quad
T^x \leq T^{y} \quad \mbox{w.p.1}.
\nonumber
\end{eqnarray}

\begin{figure}
\begin{center}
\includegraphics[width=0.5\linewidth]{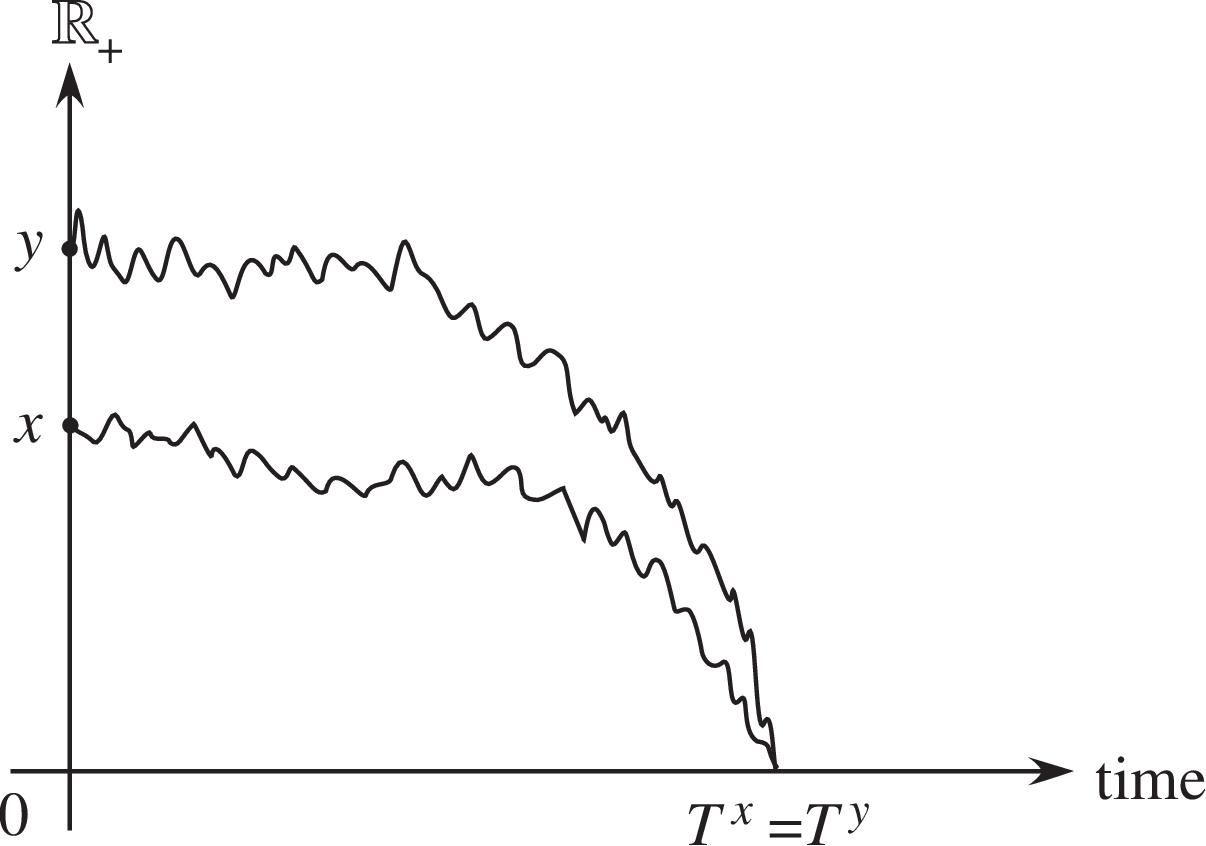}
\end{center}
\caption{\small
In the intermediate fractional dimensions,
$3/2 < D < 2$, there is a positive probability
that two Bessel processes starting from different
initial positions, $0<x<y< \infty$, return
to the origin simultaneously, $T^x=T^y$.}
\label{fig:BESflow}
\end{figure}

The interesting fact is that in the 
intermediate fractional dimensions,
$\overline{D}_{\rm c} < D < D_{\rm c}$,
it is possible to see the coincide
$T^{x}=T^{y}$ even for $x < y$.
See Fig.\ref{fig:BESflow}.

\begin{thm}
\label{thm:BESD2}
For $0 < x < y < \infty$, 
\begin{description}
\item{\rm (i)} \quad
$D \leq 3/2 \quad \Longrightarrow \quad
T^{x} < T^{y} \quad \mbox{w.p.1}$.
\item{\rm (ii)} \quad
$3/2 < D < 2 \quad \Longrightarrow \quad
\rP(T^x=T^y) > 0$.
\end{description}
\end{thm}
\vskip 0.5cm

Theorem \ref{thm:BESD2} (ii) is obtained by proving 
that, for $0 < x < y < \infty, 3/2 < D < 2$, the event
\begin{equation}
\sup_{t < T^x} \frac{X^{y}(t)-X^{x}(t)}{X^{x}(t)} < \infty
\label{eqn:sup1}
\end{equation}
occurs with positive probability.
If (\ref{eqn:sup1}) holds,
$^{\exists} c < \infty$ s.t. 
$X^{y}(t) \leq (1+c) X^{x}(t), 0 < t < T^x$
and thus 
$X^x(t)=0 \Longrightarrow X^y(t)=0
\Longrightarrow T^x=T^y$.
We can confirm that
the difference between the event $\{T^x=T^y \}$
and the event (\ref{eqn:sup1}) has probability zero
(see Section 1.10 in \cite{Law05}).

A striking fact is the following exact formula:
for $3/2 < D < 2, 0 < x < y < \infty$,
\begin{eqnarray}
\rP(T^x=T^y) &=& 1-\frac{\Gamma(D-1)}{\Gamma(2(D-1)) \Gamma(2-D)}
\left(\frac{y-x}{y} \right)^{2D-3}
\nonumber\\
&& \quad \times 
F \left( 2D-3, D-1, 2(D-1); \frac{y-x}{y} \right),
\label{eqn:Cardy1}
\end{eqnarray}
where $F(\alpha, \beta, \gamma; z)$ is {\it Gauss' hypergeometric
function}
$$
F(\alpha, \beta, \gamma; z)
=\sum_{i=0}^{\infty} \frac{(\alpha)_i (\beta)_i}
{(\gamma)_i} \frac{z^i}{i!}
$$
with the Pochhammer symbol $(c)_i=c(c+1) \cdots (c+i-1)$.
If we set $D=5/3$, (\ref{eqn:Cardy1}) becomes 
a version of exact formula for a physical quantity
called `crossing probability' 
that Cardy derived in a critical percolation model
\cite{Car92,Car01}.
Cardy's formula has been extended in the context
of SLE (see Section 6.7 of \cite{Law05}),
but I think that this exact formula
for the Bessel flow can be also called {\it Cardy's formula}.

\subsection{Schramm-Loewner evolution (SLE)
as complexification of Bessel flow}

\begin{figure}
\begin{center}
\includegraphics[width=0.5\linewidth]{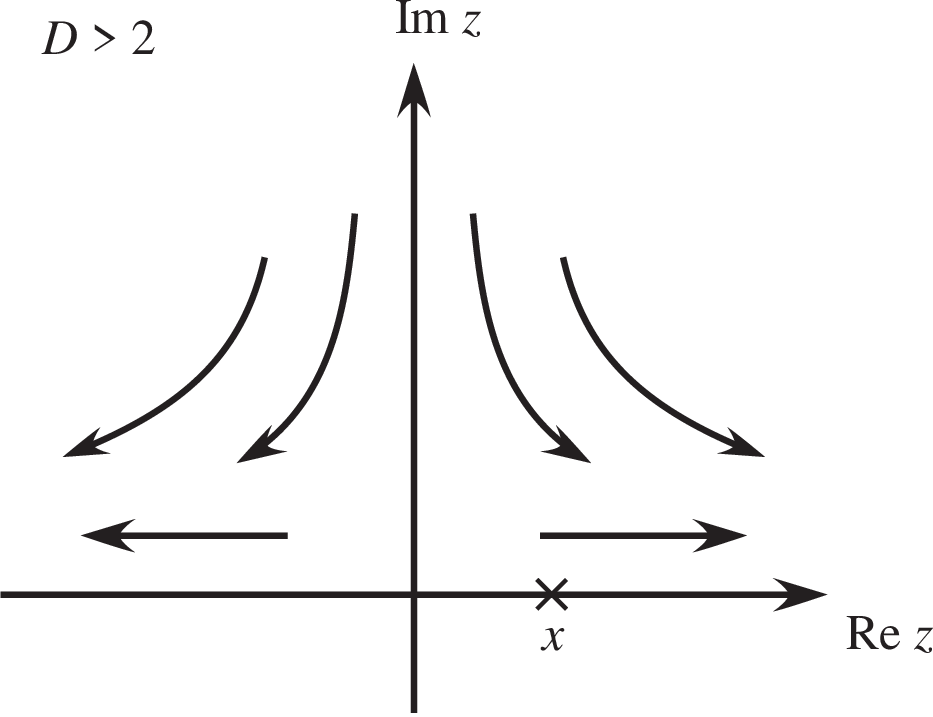}
\end{center}
\caption{\small
A schematic picture of `SLE flow' 
on $\overline{\H} \setminus \{0\}$ for $D > 2$.}
\label{fig:SLEflow}
\end{figure}

Now we consider an extension of the Bessel flow
$X^{x}(t)$ defined on $\R_+$ to flow
on the upper-half complex-plane
$\H=\{z=x+\sqrt{-1} y : x \in \R, y > 0\}$
and its boundary $\partial \H=\R$.
We set 
$Z^{z}(t)=X^{z}(t)+\sqrt{-1} Y^{z}(t) \in 
\overline{\H} \setminus \{0\}
=\H \cup \R \setminus \{0\}, t \geq 0$
and complexificate (\ref{eqn:BES2}) as
\begin{equation}
dZ^z(t)=dB(t)+ \frac{D-1}{2} \frac{dt}{Z^{z}(t)},
\quad t \in [0, \infty)
\label{eqn:SLE1}
\end{equation}
with the initial condition
$$
Z^{z}(0)=z=x+\sqrt{-1} y \in \overline{\H} \setminus \{0\}.
$$
The crucial point of this complexification of
the Bessel flow is that the BM remains real.
Then, there is asymmetry between the real part
and the imaginary part of the flow in $\H$,
\begin{eqnarray}
\label{eqn:SLE3a}
&& dX^{z}(t)=dB(t)+\frac{D-1}{2}
\frac{X^{z}(t)}{(X^z(t))^2+(Y^z(t))^2} dt, \\
\label{eqn:SLE3b}
&& dY^{z}(t)=-\frac{D-1}{2}
\frac{Y^{z}(t)}{(X^z(t))^2+(Y^z(t))^2} dt.
\end{eqnarray}
Assume $D >1$.
Then as indicated by the minus sign in the RHS of (\ref{eqn:SLE3a}),
the flow is downward in $\overline{\H}$.
If the flow goes down and arrives at the real axis,
the imaginary part vanishes, $Y^{z}(t)=0$,
and Eq.(\ref{eqn:SLE3a}) is reduced to be
the same equation as Eq.(\ref{eqn:BES2})
for the BES$^{(D)}$, which is now
considered for $\R \setminus \{0\}$.
If $D > 2$, by Theorem \ref{thm:BESD1} (ii), the flow
on $\R \setminus \{0\}$ is asymptotically outward,
$X^{x} \to \pm \infty$ as $t \to \infty$.
Therefore, the flow on $\overline{\H}$ will be described
as shown by Fig.\ref{fig:SLEflow}.
The behavior of flow should be, however, more complicated
when $\overline{D}_{\rm c}=3/2 < D < D_{\rm c}$
and $1 < D < \overline{D}_{\rm c}$.

For $z \in \overline{\H} \setminus \{0\}, t \geq 0$, let
\footnote{
Since $B^x(t), X^x(t)$ and $Z^{z}(t)$ are stochastic
processes, they are considered as functions of
time $t \geq 0$, where the initial
values $x$ and $z$ are put as superscripts
($B(t) \equiv B^{0}(t)$).
On the other hand, as explained below,
$g_t$ is considered as a conformal transformation
from a domain $H_t \subset \H$ to $\H$,
and thus it is described as a function of
$z \in H_t$; $g_t(z)$,
where time $t$ is a parameter and put as a subscript.
}
\begin{equation}
g_t(z)=Z^{z}(t)-B(t).
\label{eqn:SLE4}
\end{equation}
Then, Eq.(\ref{eqn:SLE1}) is written as follows:
\begin{equation}
\frac{d}{dt} g_t(z)=\frac{D-1}{2}
\frac{1}{g_t(z)+B(t)}, 
\quad t \geq 0,
\label{eqn:SLE5}
\end{equation}
with the initial condition
$g_0(z)=z \in \overline{\H} \setminus \{0\}$.
For each $z \in \overline{\H} \setminus \{0\}$, set
\begin{equation}
T^z= \inf \{t > 0: g_t(z)+B(t)=0\},
\label{eqn:Tz}
\end{equation}
then the solution of Eq.(\ref{eqn:SLE5}) exists up to
time $T^z$. For $t \geq 0$ we put
\begin{equation}
H_t=\{z \in \H : T^z > t\}.
\label{eqn:Ht1}
\end{equation}

This ordinary differential equation (\ref{eqn:SLE5})
involving the BM is nothing but the celebrated
{\it Schramm-Loewner evolution} (SLE) \cite{Sch00,Law05}.
For each $t \geq 0$, the solution $g_t(z)$ of (\ref{eqn:SLE5})
gives a unique conformal transformation from $H_t$ to $\H$:
$$
g_t(z) : \quad
H_t \to \H, \quad
\mbox{conformal}
$$
such that
$$
g_t(z)=z+\frac{a(t)}{z}+{\cal O} \left( \frac{1}{|z|^2} \right),
\quad z \to \infty
$$
with
$$
a(t)=\frac{D-1}{2} t.
$$

Note that the inverse map $g_t^{-1}$ from $\H$ to $H_t,
t \geq 0$, is also conformal. For each $t \geq 0$,
there exists a limit
\begin{equation}
\gamma(t)=\lim_{z \to 0, z \in \H}
g_t^{-1}(z-B(t)),
\label{eqn:gamma1}
\end{equation}
and using basic properties of BM, we can prove that
$\gamma=\gamma[0,\infty) \equiv \{\gamma(t): t \in [0, \infty)\}
\in \overline{\H}$ is a continuous path w.p.1
running from $\gamma(0)=0$ to $\gamma(\infty)=\infty$ \cite{RS05}.
The path $\gamma$ obtained from the SLE
with the parameter $D >1$ is called 
the {\it SLE$^{(D)}$ path}
\footnote{
Usual parameters used for the SLE are
$\kappa=4/(D-1)$ \cite{Sch00} or
$a=(D-1)/2=2/\kappa$ \cite{Law05}.
If we set $\widehat{g}_t(z)=\sqrt{\kappa} g_t(z)$
and $B(t)=-W(t)$ in (\ref{eqn:SLE5}),
we have the equation in the form \cite{Sch00},
$$
\frac{d}{dt} \widehat{g}_t(z)=\frac{2}
{\widehat{g}_t(z)-\sqrt{\kappa} W(t)}.
$$
Note that $\sqrt{\kappa} W(t)$ has the
same distribution with $B(\kappa t)$,
which is a time change $t \mapsto \kappa t$
of the one-dimensional standard BM.
In SLE, the Loewner chain for $\H$
is driven by a one-dimensional BM,
which is speeded up (or slowed down)
by factor $\kappa$ compared with the standard one.
(The parameter $\kappa$ is regarded as
the diffusion constant.)
}.

The dependence on $D$ of the Bessel flow
given by Theorems \ref{thm:BESD1} and \ref{thm:BESD2}
is mapped to the feature of the SLE$^{(D)}$ paths
such that there are three phases.

\begin{description}
\item{\bf [phase 1]} \,
When $D \geq D_{\rm c}=2$, the SLE$^{(D)}$ path is a 
{\it simple curve}, {\it i.e.}, $\gamma(s) \not= \gamma(t)$
for any $0 \leq s \not= t < \infty$,
and $\gamma(0, \infty) \in \H$
({\it i.e.}, $\gamma(0, \infty) \cap \R=\emptyset$).
In this phase,
$$
H_t=\H \setminus \gamma(0, t],
\quad t \geq 0.
$$
For each $t \geq 0$, $g_t$ gives a map,
which conformally erases a simple curve
$\gamma(0,t]$ from $\H$, and the image of the
tip $\gamma(t)$ of the SLE path is
$-B(t) \in \R = \partial \H$
as given by (\ref{eqn:gamma1}).
As shown by Fig.\ref{fig:SLE_P1}, it implies that the 
`SLE flow' in $\overline{\H}$ is downward
in the vertical (imaginary-axis) direction
and outward from the position $-B(t)$
in the horizontal (real-axis) direction.
Since $Z^z(t)=g_t(z)+B(t)$ by 
(\ref{eqn:SLE4}),
if we shift this figure by $B(t)$, we will have
the similar picture to Fig.\ref{fig:SLEflow} 
for the complexificated version of Bessel 
flow for $D > 2$.
\begin{figure}
\begin{center}
\includegraphics[width=0.8\linewidth]{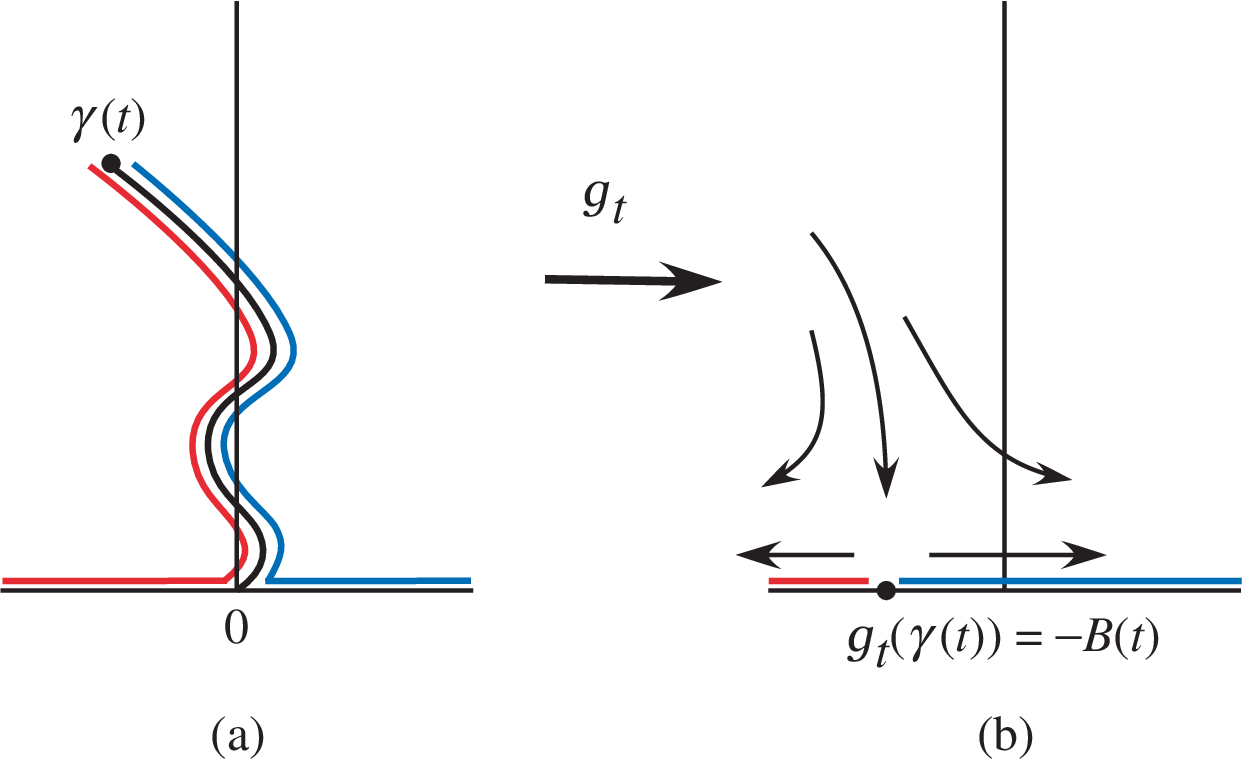}
\end{center}
\caption{\small
(a) When $D \geq 2$, the SLE$^{(D)}$ path is simple.
(b) By $g_t$, the SLE$^{(D)}$ path is erased
from $\H$. The tip of the SLE$^{(D)}$ path $\gamma(t)$
is mapped to $g_t(\gamma(t))=-B(t) \in \R$.
The flow associated by thus conformal transformation
is shown by arrows.}
\label{fig:SLE_P1}
\end{figure}

\item{\bf [phase 2]} \,
When $\overline{D}_{\rm c}=3/2 < D < D_{\rm c}=2$,
the SLE$^{(D)}$ path can osculate the real axis,
$\rP(\gamma(0,t] \cap \R \not= \emptyset) >0,
\forall t >0$.
Fig.\ref{fig:SLE_P2} (a) illustrates the moment $t>0$
such that the tip of SLE$^{(D)}$ path just 
osculates the real axis.
The closed region encircled by the path
$\gamma(0,t)$ and the line $[\gamma(t), 0] \in \R$
is called an {\it SLE hull} at time $t$ and 
denoted by $K_t$.
In this phase
$$
H_t=\H \setminus K_t, \quad
t \geq 0.
$$
That is, $g_t(z)$ is a map which erases conformally
the SLE hull from $\H$.
We can think that by this transformation
all the points in $K_t$ are simultaneously mapped to
a single point $-B(t) \in \R$,
which is the image of the tip $\gamma(t)$.
(We say that the hull $K_t$ is {\it swallowed}.
See Fig.\ref{fig:SLE_P2} (b).)
By definition (\ref{eqn:Ht1}), the moment
when $K_t=\H \setminus H_t$ is swallowed is the time
$T^z$ at which the equality
$Z^z(t)=g_t(z)+B(t)=0$ holds 
$^{\forall}z \in K_t$.
(Then the RHS of (\ref{eqn:SLE5}) diverges and all the points
$z \in K_t$ are lost from the domain of the map $g_t$.)
Theorem \ref{thm:BESD2} (ii) states that,
when $\overline{D}_{\rm c} < D < D_{\rm c}$, 
two BES$^{(D)}$ starting from different points $0<x<y< \infty$
can simultaneously return to the origin.
In the complexificated version, all $Z^z(t)$
starting from $z \in K_t$ can arrive at the origin
simultaneously
({\it i.e.}, they are all swallowed).

Osculation of the SLE path with $\R$ means
that the SLE path has loops.
Figure \ref{fig:SLE_P2B} (a) shows the event that the SLE path
makes a loop at time $t>0$.
The SLE hull $K_t$ consists of the closed region
encircled by the loop and the segment
of the SLE path between the origin and the
osculating point, and it is completely erased
by the conformal transformation $g_t$ to $\H$
as shown by Fig.\ref{fig:SLE_P2B}(b).
Let $0 < s < t$ and consider the map $g_s$,
which is the solution of (\ref{eqn:SLE5}) at time $s$.
Assume that $\gamma(s)$ is located on the loop part
of $\gamma[0, t]$ as shown by Fig.\ref{fig:SLE_P2B}(a).
The segment $\gamma[0,s]$ of the SLE path is mapped by
$g_s$ to a part of $\R$.
Since $\gamma(t)$ osculates a point in $\gamma[0,s]$,
its image $g_s(\gamma(t))$ should osculate the
real axis $\R$ as shown by Fig.\ref{fig:SLE_P2B}(c).
Since $g_s^{-1}$ is uniquely determined from $g_s$,
the above argument can be reversed.
Then equivalence between osculation of the SLE path
with $\R$ and self-intersection of the SLE path
is concluded.
In this intermediate phase $\overline{D}_{\rm c}
< D < D_{\rm c}$,
\begin{eqnarray}
&& \mbox{SLE$^{(D)}$ path $\gamma$ is {\it self-intersecting, and}}
\nonumber\\
&& \bigcup_{t > 0} \overline{K_t}=
\overline{\H} \quad \mbox{but} \quad
\gamma[0, \infty) \cap \H
\not= \H \quad \mbox{w.p.1}.
\nonumber
\end{eqnarray}

\begin{figure}
\begin{center}
\includegraphics[width=0.8\linewidth]{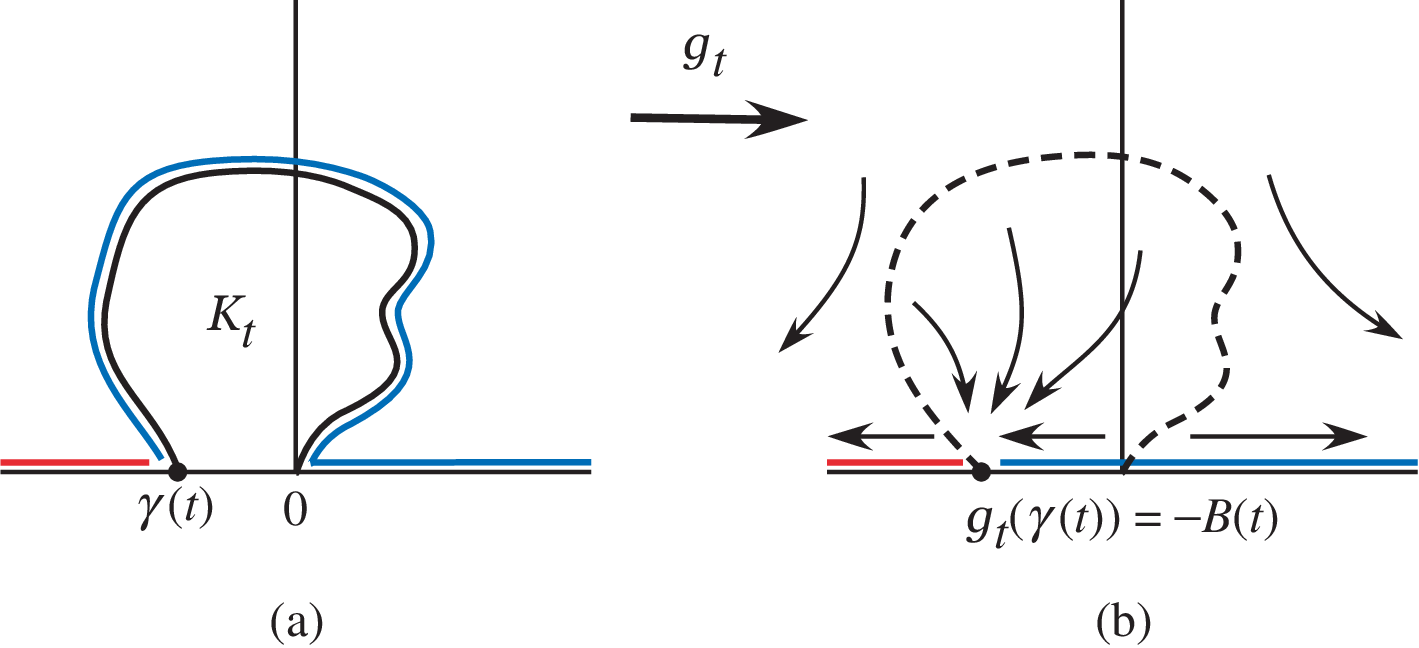}
\end{center}
\caption{\small
(a) When $3/2 < D < 2$, the SLE$^{(D)}$ path can osculate
the real axis. The SLE hull is denoted by $K_t$.
(b) The SLE hull $K_t$ is swallowed.
It means that all the points in $K_t$ are simultaneously
mapped to a single point $-B(t) \in \R$,
which is the image of the tip of the SLE$^{(D)}$ path
$\gamma(t)$.}
\label{fig:SLE_P2}
\end{figure}
\begin{figure}
\begin{center}
\includegraphics[width=0.7\linewidth]{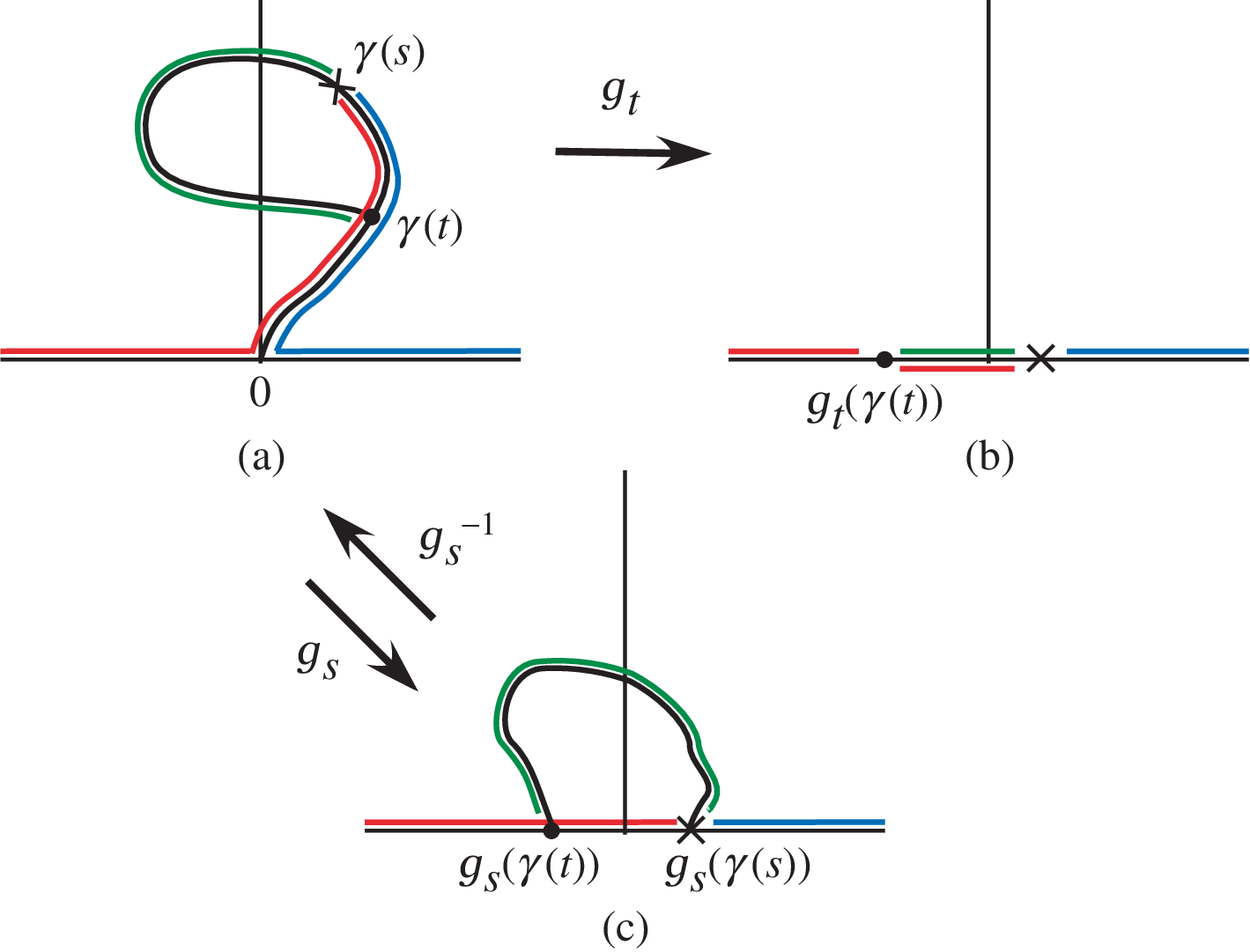}
\end{center}
\caption{\small
The event that the SLE$^{(D)}$ path osculates $\R$
is equivalent with the event that the SLE$^{(D)}$ path
makes a loop.}
\label{fig:SLE_P2B}
\end{figure}

\item{\bf [phase 3]} \,
When $1 < D \leq \overline{D}_{\rm c}=3/2$,
Theorem \ref{thm:BESD2} (i) states for the
Bessel flow that the ordering
$T^x< T^y$ is conserved for any $0 < x < y$.
It implies that in this phase
the SLE path should be a {\it space-filling curve};
$$
\gamma[0, \infty) = \overline{\H}.
$$
(Otherwise, swallow of regions occurs,
contradicting Theorem \ref{thm:BESD2} (i).)
\end{description}

\begin{figure}[htbp]
  \begin{center}
   \includegraphics[width=0.9\linewidth]{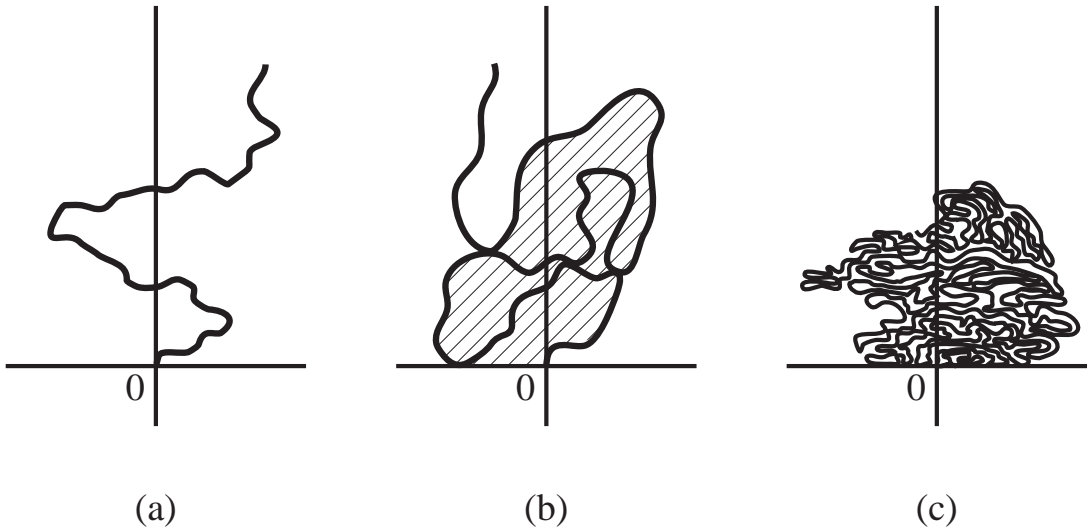}
  \end{center}
  \caption{\small
  Schematic pictures of SLE$^{(D)}$ paths in
  (a) phase 1 ($D \geq D_{\rm c}=2$),
  (b) phase 2 ($\overline{D}_{\rm c}=3/2 < D < D_{\rm c}=2$), and
  (c) phase 3 ($1 < D \leq \overline{D}_{\rm c}=3/2$).
  }
  \label{fig:SLEphases}
\end{figure}

\noindent
Figure \ref{fig:SLEphases}
summarizes the three phases of SLE paths.

By complexification of Bessel flow,
we can discuss flows on a two-dimensional plane $\overline{\H}$.
By this procedure random curves (the SLE paths) 
are generated in the plane.
The SLE paths are fractal curves and their {\it Hausdorff dimensions}
$d_{\rm H}^{(D)}$ are determined by Beffara \cite{Bef08}.
We note that a reciprocity relation is found
between $D$ and $d_{\rm H}^{(D)}$;
\begin{equation}
(D-1)(d_{\rm H}^{(D)}-1)=\frac{1}{2},
\quad D \geq \overline{D}_{\rm c}=\frac{3}{2}.
\label{eqn:Hausdorff}
\end{equation}
(In the phase 3, $D \leq \overline{D}_{\rm c}=3/2$,
$d_{\rm H}^{(D)} \equiv 2$.)

Remark that the SLE map $\gamma_t$ as well as
the SLE path $\gamma$ are functionals of the BM.
Therefore, we have statistical ensembles
of random curves $\{\gamma\}$ in the probability
space $(\Omega, {\cal F}, \rP)$.
The important consequence from the facts that
the BM is a strong Markov process with 
independent increments and $g_t$ gives
a conformal transformation is that
the statistics of $\{\gamma\}$
has a kind of stationary Markov property
(called the {\it domain Markov property}) and 
{\it conformal invariance} with respect to
transformation of the domain in which
the SLE path $\gamma$ is defined \cite{Law05}.

The highlight of the theory of SLE would be
that, if the value $D$ is properly chosen,
the statistics of $\{\gamma\}$ realizes
that of the scaling limit of important
{\it statistical mechanics model} exhibiting
{\it critical phenomena} and
fractal structures defined on
an infinite discrete lattice.
The following is a list of the correspondence
between the SLE$^{(D)}$ paths with
specified values of $D$ and
the scaling limits of models studied
in statistical mechanics and fractal physics
\footnote{
SLE$^{(D)}$ has a special property called
the {\it restriction property} iff $D=5/2 \, (\kappa=8/3)$.
It is well-know that the self-avoiding walk (SAW) model,
which has been studied as a model for polymers, 
has this property. The conformal invariance
of the scaling limit of SAW is, however, 
not yet proved. If it is proved, 
the equivalence in probability law 
between the scaling limit of SAW
and the SLE$^{(5/2)}$ path will be concluded.
}.
\begin{eqnarray}
\mbox{SLE$^{(3/2)}$}
&\Longleftrightarrow& 
\mbox{uniform spanning tree} \, \cite{LSW04}
\nonumber\\
\mbox{SLE$^{(5/3)}$}
&\Longleftrightarrow& 
\mbox{critical percolation model}
\nonumber\\
&& 
\mbox{(percolation exploration process)} \, \cite{Smi01}
\nonumber\\
\mbox{SLE$^{(2)}$}
&\Longleftrightarrow& 
\mbox{Gaussian free surface model (contour line)} \, \cite{SS05}
\nonumber\\
\mbox{SLE$^{(7/3)}$}
&\Longleftrightarrow& 
\mbox{critical Ising model (Ising interface)} \, \cite{Smi10}
\nonumber\\
\mbox{SLE$^{(5/2)}$}
&\Longleftrightarrow& 
\mbox{self-avoiding walk [conjecture]}
\nonumber\\
\mbox{SLE$^{(3)}$}
&\Longleftrightarrow& 
\mbox{loop-erased random walk} \, \cite{LSW04}
\nonumber
\end{eqnarray}

\subsection{Dyson's BM model
as multivariate extension of Bessel process}

Here we consider stochastic motion of two particles
$(X_1(t), X_2(t))$ in one dimension $\R$
satisfying the following SDEs,
\begin{eqnarray}
&& dX_1(t) = dB_1(t)
+\frac{\beta}{2} \frac{dt}{X_1(t)-X_2(t)}
\nonumber\\
&& dX_2(t) = dB_2(t)
+\frac{\beta}{2} \frac{dt}{X_2(t)-X_1(t)}
\quad t \geq 0,
\label{eqn:Dyson1}
\end{eqnarray}
with the initial condition $X_1(0) < X_2(0)$, 
where $B_1(t)$ and $B_2(t)$ are independent
one-dimensional standard BMs and $\beta >0$
is a `coupling constant' of the two particles.
The second terms in (\ref{eqn:Dyson1})
represent the repulsive force acting between two 
particles, which is proportional to the inverse
of distance $X_2(t)-X_1(t)$ of the two particles.
Since it is a central force ({\it i.e.}, 
depending only on distance, and thus symmetric
for two particles), the `center of mass'
$X_{\rm c}(t) \equiv (X_2(t)+X_1(t))/2$
is proportional to a BM;
$X_{\rm c}(t) \d= B(t)/\sqrt{2} \d= B(t/2), ^{\forall} t \geq 0$,
where $B(t)$ is a one-dimensional
standard BM different from $B_1(t)$ and $B_2(t)$
and the symbol $\d=$ denotes the equivalence in distribution.
(Note that variance (quadratic variation) $(dX_{\rm c}(t))^2=
\{(dB_1(t)+dB_2(t))/2\}^2=dt/2$,
since $(dB_1(t))^2=(dB_2(t))^2=dt$ and
$dB_1(t) dB_2(t)=0$.)
On the other hand, we can see that the relative
coordinate defined by
$X_{\rm r}(t) \equiv (X_2(t)-X_1(t))/\sqrt{2}$
satisfies the SDE
$$
dX_{\rm r}(t)=d \widetilde{B}(t)
+\frac{\beta}{2} \frac{dt}{X_{\rm r}(t)},
\quad t \geq 0,
$$
where $\widetilde{B}(t)$ is a BM different from
$B_1(t), B_2(t), B(t)$.
It is nothing but the SDE for BES$^{(D)}$ with
$D=\beta+1$.

Dyson \cite{Dys62} introduced $N$-particle systems
of interacting BMs in $\R$  
as a solution 
$\X(t)=(X_1(t), X_2(t), \dots, X_N(t))$
of the following system of SDEs,
\begin{equation}
dX_i(t)=dB_i(t)+ \frac{\beta}{2}
\sum_{1 \leq j \leq N: j \not=i}
\frac{dt}{X_i(t)-X_j(t)},
\quad t \in (0, \tau),
\quad i=1,2, \dots, N, 
\label{eqn:Dyson3}
\end{equation}
where $\{B_i(t) \}_{i=1}^N$ are independent
one-dimensional standard BMs
and we define
\begin{eqnarray}
&& \sigma_{ij}=\inf \{ t \geq 0 :
X_i(t) \not= X_j(t)\}, \quad 1 \leq i < j \leq N, \nonumber\\
&& \tau_{ij}=\inf \{ t > \sigma_{ij} :
X_i(t)=X_j(t)\}, \quad 1 \leq i < j \leq N,
\nonumber\\
&& \qquad \qquad \tau=\min_{1 \leq i < j \leq N} \tau_{ij}.
\nonumber
\end{eqnarray}
It is called {\it Dyson's BM model} with parameter $\beta$
\cite{Meh04,For10}.
As shown above, the $N=2$ case of Dyson's BM model
is a coordinate transformation of the pair of
a (time-change of) BM and BES$^{(\beta+1)}$.
In this sense, Dyson's BM model can be regarded
as a multivariate (multi-dimensional) extension
of BES$^{(\beta+1)}, \beta>0$
\footnote{
We can prove that $\tau < \infty$ for $\beta <1$
and $\tau = \infty$ for $\beta \geq 1$ \cite{RS93}.
The critical value $\beta_{\rm c}=1$
corresponds to $D_{\rm c}=2$ of
BES$^{(D)}$.
}.
In particular, we will characterize Dyson's BM model
with $\beta=2$ as an extension of the three-dimensional
Bessel process, BES$^{(3)}$, in Section 2.

\SSC{Two aspects of the Dyson model}

In this section, we study the special case of
Dyson's BM model with parameter $\beta=2$.
We call this special case simply {\it the Dyson model}
\cite{KT10b}.
As shown above, the case $\beta=2$ corresponds to the case $D=3$
of Bessel process. In Sect.1.2, we have shown that
BES$^{(3)}$ has two aspects;
{\bf (Aspect 1)} as a radial coordinate of three-dimensional
BM, and
{\bf (Aspect 2)} as a one-dimensional BM conditioned
to stay positive.
We show that the Dyson model inherits these two aspects
from BES$^{(3)}$ \cite{KT_Sugaku}.

\subsection{The Dyson model as eigenvalue process}

Dyson introduced the process (\ref{eqn:Dyson3})
with $\beta=1, 2$, and 4 as the eigenvalue processes
of matrix-valued BMs in the Gaussian orthogonal 
ensemble (GOE), the Gaussian unitary ensemble (GUE),
and the Gaussian symplectic ensemble (GSE) \cite{Dys62,Meh04,For10}
\footnote{
Precisely speaking, Dyson considered the Ornstein-Uhlenbeck
processes of eigenvalues such that
as stationary states they have the eigenvalue distributions
of random matrices in GOE, GUE, and GSE.
Here we consider matrix-valued BMs, so 
variances increase in proportion to time $t \geq 0$.
}.

For $\beta=2$ with given $N \in \N$, we prepare
$N$-tuples of one-dimensional standard BMs
$\{B_{ii}^{x_i} \}_{j=1}^{N}$,
each of which starts from $x_i \in \R$,
$N(N-1)/2$-tuples of pairs of BMs
$\{B_{ij}, \widetilde{B}_{ij}\}_{1 \leq i< j \leq N}$,
all of which start from the origin, where
the totally $N+2 \times N(N-1)/2=N^2$ BMs
are independent from each other.
Then consider an $N \times N$ Hermitian-matrix-valued BM
\footnote{
In usual Gaussian random matrix ensembles,
mean is assumed to be zero.
The corresponding matrix-valued BM are then considered
to be started from a zero matrix,
{\it i.e.}, $x_i=0, 1 \leq i \leq N$
in (\ref{eqn:Mx0}).
In random matrix theory, general case with
non-zero means ({\it i.e.}, $x_i \not= 0$) 
is discussed with the terminology
`random matrices in an external source'
(see, for example, \cite{BH98}).
From the view point of stochastic processes,
imposing external sources 
to break symmetry of the system
corresponds to changing initial state.
}.
\begin{equation}
\cH^{\x}(t) = \left( \begin{array}{cccc}
B_{11}^{x_1}(t) 
& \displaystyle{\frac{B_{12}(t)+\sqrt{-1} \widetilde{B}_{12}(t)}{\sqrt{2}}}
& \cdots & 
\displaystyle{\frac{B_{1N}(t)+\sqrt{-1} \widetilde{B}_{1N}(t)}{\sqrt{2}}}
\cr
\displaystyle{\frac{B_{12}(t)-\sqrt{-1} \widetilde{B}_{12}(t)}{\sqrt{2}}}
& B_{22}^{x_2}(t) & \cdots &
\displaystyle{\frac{B_{2N}(t)+\sqrt{-1} \widetilde{B}_{2N}(t)}{\sqrt{2}}}
\cr
\cdots & \cdots & \cdots & \cdots \cr
\displaystyle{\frac{B_{1N}(t)-\sqrt{-1} \widetilde{B}_{1N}(t)}{\sqrt{2}}}
& 
\displaystyle{\frac{B_{2N}(t)-\sqrt{-1} \widetilde{B}_{2N}(t)}{\sqrt{2}}}
 & \cdots & B_{NN}^{x_N}(t) 
\end{array} \right).
\label{eqn:Mxt}
\end{equation}
By this definition, the initial state of this BM is
given by the diagonal matrix 
\begin{equation}
\cH^{\x}(0)={\rm diag} (x_1, x_2, \dots, x_N).
\label{eqn:Mx0}
\end{equation}
We assume $x_1 \leq x_2 \leq \cdots \leq x_N$.

Remember that when we introduced BES$^{(D)}$ in Sect.1.2,
we considered the $D$-dimensional vector-valued BM,
(\ref{eqn:DBM}), by preparing $D$-tuples of
independent one-dimensional standard BMs for its elements.
Since the dimension of the space ${\rm H}(N)$
of $N \times N$ Hermitian matrices is
${\rm dim} \, {\rm H}(N)=N^2$,
we need $N^2$ independent BMs for elements
to describe a BM in this space ${\rm H}(N)$.

Corresponding to calculating an absolute value of 
$\B^{\x}(t)$, by which BES$^{(D)}$ was introduced as
(\ref{eqn:BES1}), here we calculate eigenvalues of
$\cH^{\x}(t)$.
For any $t \geq 0$, there is a family of $N \times N$
unitary matrices $\{U(t)\}$ which diagonalize $\cH^{\x}(t)$,
$$
U(t)^{*} \cH^{\x}(t) U(t) =
{\rm diag} (\lambda_1(t), \dots, \lambda_N(t)),
\quad t \geq 0.
$$
Let $\W_N^{\rm A}$ be the
Weyl chamber of type A$_{N-1}$ defined by
$$
\W_N^{\rm A}
\equiv \{\x=(x_1, x_2, \dots, x_N)
\in \R^N: x_1 < x_2 < \cdots < x_N\}.
$$
If we impose the condition
$\vlambda(t) \equiv (\lambda_1(t), \dots, \lambda_N(t)) \in
\W_N^{\rm A}, t \geq 0$,
$U(t)$ is uniquely determined.

For each $t \geq 0$, $\lambda_i(t), 1 \leq i \leq N$
are functionals of $\{B_{ii}^{x_i}(t),
B_{ij}(t), \widetilde{B}_{ij} \}_{1 \leq i < j \leq N}$,
and then, again we can apply It\^o's formula
to take into account propagation of error
correctly to derive SDEs for the eigenvalue process
$\vlambda(t), t \geq 0$ \cite{Bru89,Bru91,KT04}.
The result is the following,
\begin{equation}
d \lambda_i(t)=dB_i^{x_i}(t)
+ \sum_{1 \leq j \leq N: j \not=i}
\frac{dt}{\lambda_i(t)-\lambda_j(t)},
\quad t \in (0, \infty), \quad i=1,2, \dots, N,
\label{eqn:lambda}
\end{equation}
where $\{B_i^{x_i} \}_{i=1}^{N}$ are independent
BM different from the BMs used to define $\cH^{\x}(t)$
by (\ref{eqn:Mxt}).
It is indeed the $\beta=2$ case of (\ref{eqn:Dyson3})
as derived by Dyson 
(originally not by such a stochastic calculus
but by applying the perturbation theory in 
quantum mechanics) \cite{Dys62}.

Now the correspondence is summarized as follows.
$$
\begin{array}{lccl}
{\mbox{\bf [Aspect 1]}} &  &  &  \cr
  & \mbox{BES$^{(3)}$} & \Longleftrightarrow & 
  \mbox{radial coordinate of} \cr
  &                    &                     &
  \mbox{$D=3$ vector-valued BM} \cr
  &  &  & \cr
  & \mbox{the Dyson model} & \Longleftrightarrow &
  \mbox{eigenvalue process of} \cr
  & \mbox{with $N$ particles} & &
  \mbox{$N \times N$ Hermitian-matrix-valued BM}
\end{array}
$$

\subsection{The Dyson model as noncolliding BM}

Here we try to extend the formula (\ref{eqn:p3})
for BES$^{(3)}$ to multivariate versions.

First we consider a set of two operations,
identity ($\sigma_1=$ id) 
and reflection ($\sigma_2=$ ref), such that
for $x \in \R$, $\sigma_1(x)=x$ and
$\sigma_2(x)=-x$, and signatures are given as
${\rm sgn}(\sigma_1)=1$ and ${\rm sgn}(\sigma_2)=-1$,
respectively.
Then we have
\begin{equation}
p_t^{\rm abs}(y|x)
=\sum_{\sigma \in \{{\rm id}, {\rm ref}\}}
{\rm sgn}(\sigma) p_{t}(y|\sigma(x)),
\quad t \geq 0, \quad x, y, \in \R.
\label{eqn:p3b}
\end{equation}
Next we consider a set of all permutations of $N$ indices
$\{1,2, \dots, N\}$, which is denoted by ${\cal S}_N$,
and put the following multivariate function
\begin{equation}
\sum_{\sigma \in {\cal S}_N}
{\rm sgn}(\sigma) \prod_{i=1}^{N} 
p_t(y_{\sigma(i)}|x_i)
=\det_{1 \leq i, j \leq N} 
\Big[ p_t(y_i|x_j) \Big]
\label{eqn:KM1}
\end{equation}
of $\x=(x_1, \dots, x_N) \in \W_N^{\rm A}$ and
$\x=(y_1, \dots, y_N) \in \W_N^{\rm A}$ with a parameter $t \geq 0$.
Following the argument by Karlin and McGregor \cite{KM59},
we can prove that this determinant gives the 
transition probability density with duration $t$
from the state $\x$ to the state $\y$ of $N$-dimensional
absorbing BM, 
$\B^{\x}(t)=(B_1^{x_1}(t), \dots, B_N^{x_N}(t))$, 
in a domain $\W_N^{\rm A}$, in which absorbing walls
are put at the boundaries of $\W_N^{\rm A}$.
Since the boundaries of $\W_N^{\rm A}$ are the hyperplanes
$x_i=x_j, 1 \leq i < j \leq N$,
the Brownian particle $\B^{\x}(t)$ is annihilated,
when any coincidence of the values of coordinates
of $\B^{\x}(t)$ occurs.
The `survival probability' that the BM is not yet
absorbed at the boundary is a monotonically decreasing
function of time.
If we regard the $i$-th coordinate $B_i^{x_i}(t)$
as the position of $i$-th particle on $\R, 1 \leq i \leq N$,
the state $\x \in \W_N^{\rm A}$ is considered 
to represent a configuration of $N$ particles
on $\R$ such that a strict ordering $x_1 < x_2 < \cdots < x_N$
of positions is maintained,
while the state absorbed at the boundary, 
$\x \in \partial \W_M^{\rm A}$, is a configuration
in which collision occurs between some
pair of neighboring pairs of particles;
$1 \leq ^{\exists} i \leq N-1$,
s.t. $x_i=x_{i+1}$.
Since if such collision occurs, the process 
is totally annihilated, this many-particle
system is called {\it vicious walker model} 
\cite{Fis84,KGV00,Joh02,KT02,Joh03,CK03}.

As already noted below Eq.(\ref{eqn:p3}),
$h_1(x)=x$ is a harmonic function in $(0, \infty)$ 
conditioned $\phi(0)=0$.
Similarly, if we consider a harmonic function of
$N$ variables
\begin{equation}
\Delta^{(N)} h_N(\x) 
\equiv \sum_{i=1}^{N} \frac{\partial^2}{\partial x_i^2}
h_N(\x)=0
\label{eqn:hN1}
\end{equation}
conditioned $h_N(\widehat{\x})=0, \,
^{\forall} \widehat{\x} \in \partial \W_N^{\rm A}$,
we will have the following {\it product of differences}
\begin{equation}
h_N(\x)=\prod_{1 \leq i < j \leq N}
(x_j-x_i),
\label{eqn:hN2}
\end{equation}
which is identified with the
{\it Vandermonde determinant} 
$\det_{1 \leq i, j \leq N}[x_i^{j-1}]$.

Combining above consideration, we put the following
function
\begin{equation}
p_t^{N}(\y|\x)=\frac{h_N(\y)}{h_N(\x)}
\det_{1 \leq i, j \leq N} 
\Big[ p_t(y_i|x_j) \Big],
\quad t \geq 0, \quad \x, \y \in \W_N^{\rm A}.
\label{eqn:nonc}
\end{equation}
We can show that the factor $h_N(\y)/h_N(\x)$
provides an exact renormalization of the vicious walker model
(the absorbing BM in $\W_N^{\rm A}$) to compensate any decay of 
total mass of the process by collision 
(by absorption at $\partial \W_N^{\rm A}$), and that
(\ref{eqn:nonc}) gives the transition probability density
function for the $N$-particle system of one-dimensional
BMs {\it conditioned never to collide with each other
forever} (which we simply call the
{\it noncolliding BM}) \cite{Gra99,KT02,KT03a}.

Moreover, by using the harmonicity (\ref{eqn:hN1}),
we can confirm that (\ref{eqn:nonc}) satisfies the
following partial differential equation,
\begin{equation}
\frac{\partial}{\partial t}
p_t^N(\y|\x)
=\frac{1}{2} \Delta^{(N)} p_t^N(\y|\x)
+\sum_{1 \leq i \not= j \leq N}
\frac{1}{x_i-x_j} \frac{\partial}{\partial x_i}
p_t^N(\y|\x)
\label{eqn:nonc2}
\end{equation}
with the initial condition $p_0^N(\y|\x)=\delta(\y-\x)
=\prod_{i=1}^{N} \delta(y_i-x_i)$ \cite{KT02}.
It can be regarded as the backward Kolmogorov equation
of the stochastic process with $N$ particles,
$\X(t)=(X_1(t), \dots, X_N(t))$, 
which solves the system of SDEs,
\begin{equation}
dX_i(t)=dB_i^{x_i}(t)
+\sum_{1 \leq j \leq N: j \not=i}
\frac{dt}{X_i(t)-X_j(t)}, \quad
t \in (0, \infty), \quad i=1,2, \dots, N.
\label{eqn:nonc3}
\end{equation}

Eq.(\ref{eqn:nonc3}) is identified with the $\beta=2$
case of (\ref{eqn:Dyson3}).
Then the equivalence between the Dyson model
and the noncolliding BM is proved.

The result is summarized as follows.
$$
\begin{array}{lccl}
{\mbox{\bf [Aspect 2]}} &  &  &  \cr
  & \mbox{BES$^{(3)}$} & \Longleftrightarrow & 
  \mbox{$h$-transform of absorbing BM in $(0, \infty)$} \cr
  &                    & \Longleftrightarrow &
  \mbox{BM conditioned to stay positive} \cr
  &  &  & \cr
  & \mbox{the Dyson model} & \Longleftrightarrow &
  \mbox{$h$-transform of absorbing BM in $\W_N^{\rm A}$} \cr
  & & \Longleftrightarrow &
  \mbox{noncolliding BM}
\end{array}
$$

The fact that BES$^{(3)}$ and the Dyson model have two
aspects implies useful relation between 
projection from higher dimensional spaces
and restriction by imposing conditions.
For matrix-valued processes, projection is
performed by integration over irrelevant components, 
and noncolliding conditions are generally 
expressed by the Karlin-McGregor determinants.
The processes discussed here are temporally
homogeneous ones, but we can also discuss two
aspects of temporally inhomogeneous processes.
Actually, we have shown that, from the fact that
the temporally inhomogeneous version of noncolliding
BM has the two aspects, 
the Harish-Chandra-Itzykson-Zuber integral formula
\cite{HC57,IZ80} is derived,
\begin{equation}
\int_{{\rm U}(N)} dU \,
\exp \left\{ - \frac{1}{2 \sigma^2}
{\rm Tr}(\Lambda_{\x}-U^{*} \Lambda_{\y} U)^2 \right\}
=\frac{C_N \sigma^{N^2}}{h_N(\x) h_N(\y)}
\det_{1 \leq i, j \leq N}
\Big[ p_t(y_i|x_j) \Big],
\label{eqn:HC1}
\end{equation}
$\sigma^2>0$,
where $dU$ denotes the Haar measure of ${\rm U}(N)$
normalized as $\int_{{\rm U}(N)} dU=1$,
$\Lambda_{\x}={\rm diag}(x_1, \dots, x_N), 
\Lambda_{\y}={\rm diag}(y_1, \dots, y_N)$
with $\x, \y \in \W_N^{\rm A}$,
and $C_N=(2\pi)^{N/2} \prod_{i=1}^{N} \Gamma(i)$
\cite{KT03b}
\footnote{
We can apply the present argument also for
processes associated with Weyl chambers of
other types. See \cite{KTNK03,KT04}.
}.

\SSC{Determinantal processes and entire functions}
\subsection{Aspect 2 of the Dyson model}

As Aspect 2, the Dyson model is constructed as the
$h$-transform of the absorbing BM in $\W_N^{\rm A}$.
Therefore, at any positive time $t > 0$ the configuration
is given as an element of $\W_N^{\rm A}$,
\begin{equation}
\X(t)=(X_1(t), X_2(t), \dots, X_N(t)) \in
\W_N^{\rm A}, \quad t > 0,
\label{eqn:Xt}
\end{equation}
and there is no multiple point at which
coincidence of particle positions, 
$X_i(t)=X_j(t), i \not= j$, occurs.
That is, the Dyson model is equivalent with
the noncolliding BM.
We can consider the Dyson model, however, 
starting from initial configurations with multiple points.
In order to describe configurations with multiple points, 
we represent each particle configuration by a sum 
of delta measures in the form
\begin{equation}
\xi(\cdot)=\sum_{i \in \I} \delta_{x_i}(\cdot),
\quad t \geq 0
\label{eqn:xi}
\end{equation}
with a sequence of points in $\R,
\x=(x_i)_{i \in \I}$,
where $\I$ is a countable index set.
Here for $y \in \R$,
$\delta_y(\cdot)$ denotes the delta measure such that
$\delta_y(x)=1$ for $x=y$ and
$\delta_y(x)=0$ otherwise.
Then, for (\ref{eqn:xi}) and $A \subset \R$,
$\xi(A) \equiv \int_{A} \xi(dx)
=\sum_{i \in \I: x_i \in A} 1
=\sharp\{ x_i: x_i \in A\}$.
If the total number of particles $N$ is finite,
$\I=\{1,2, \dots, N\}$, but we would like to also
consider the cases with $N=\infty$.
We call measures of the form (\ref{eqn:xi})
satisfying the condition
$\xi(K) < \infty$ for any compact subset $K \subset \R$
nonnegative integer-valued Radon measures on $\R$
and write the space of them as $\mM$. 
The set of configurations without multiple point
is denoted by 
$\mM_0=\{\xi \in \mM : \xi(\{x\}) \leq 1, ^{\forall} x \in \R\}$.
There is a trivial correspondence between $\W_N^{\rm A}$
and $\mM_0$.

First we assume $\xi=\sum_{i \in \I} \delta_{x_i} \in \mM_0$,
$\xi(\R)=N \in \N$
and consider the Dyson model as an $\mM_0$-valued 
diffusion process,
\begin{equation}
\Xi(t, \cdot)=\sum_{i \in \I} \delta_{X_i(t)} (\cdot),
\quad t \geq 0, 
\label{eqn:Xit}
\end{equation}
starting from the initial configuration $\xi$,
where $\X(t)=(X_1(t), \cdots, X_N(t))$ is the
solution of (\ref{eqn:nonc3}) under the initial 
configuration $\x=(x_1, \dots, x_N) \in \W_N^{\rm A}$.
We write the process as $(\Xi(t), \P^{\xi})$
and express the expectation with respect to
the probability law $\P^{\xi}$ of the Dyson model by
$\E^{\xi}[\, \cdot \,]$.
We introduce a filtration $\{\F(t)\}_{t \in [0, \infty)}$
on the space of continuous paths
$\c([0, \infty) \to \mM)$ defined by
$\F(t)=\sigma(\Xi(s), s \in [0, t])$, where
$\sigma$ denotes the sigma field.

Then we introduce a sequence of independent
one-dimensional standard BMs,
$\B^{\x}(t)=(B_i^{x_i}(t))_{i \in \I},
t \geq 0$
and write the expectation with respect to them as 
$\rE^{\x}[\, \cdot \,]$.

Let $\1_{\W_N^{\rm A}}(\x)=1$ if $\x \in \W_N^{\rm A}$
and $\1_{\W_N^{\rm A}}(\x)=0$ otherwise.
Then Aspect 2 of the Dyson model is expressed by
the following equality;
for any $0 < t < T < \infty$,
any symmetric function $g$ on $\R^{N}$,
\begin{equation}
\E^{\xi}[g(\X(t))]
=\rE^{\x} \left[ g(\B(t)) \1_{\W_N^{\rm A}}(\B(t))
\frac{|h_N(\B(t))|}{h_N(\x)} \right],
\label{eqn:Asp2_1}
\end{equation}
where we have assumed the relations
$\xi=\sum_{i \in \I} \delta_{x_i} \in \mM_0,
\xi(\R)=N \in \N,
\x=(x_1, \dots, x_N) \in \W_N^{\rm A}$
and (\ref{eqn:Xit}).
Note that the indicator $\1_{\W_N^{\rm A}}(\B(t))$
in the RHS annihilates the BM, $\B(t)$, if it hits
any boundary of the Weyl chamber, $\partial \W_N^{\rm A}$,
and the factor $|h_N(\B(t))|/h_N(\x)$ performs
the $h$-transform of measure for Brownian paths.
That is, the RHS of (\ref{eqn:Asp2_1}) indeed
gives the expectation of $g$ with respect to the
process obtained by the $h$-transform of the absorbing
BM in $\W_N^{\rm A}$.

If we apply the Karlin-McGregor determinantal formula
(see (\ref{eqn:KM1})),
$$
\mbox{(RHS)}=
\rE^{\x} \left[ \sum_{\sigma \in \S_N}
{\rm sgn}(\sigma) \1_{\W_N^{\rm A}}(\sigma(\B(t)))
g(\B(t))
\frac{|h_N(\B(t))|}{h_N(\x)} \right],
$$
where we recall the definition of determinant
and let $\sigma(\B(t))=(B_{\sigma(1)}, \dots,
B_{\sigma(N)}), \sigma \in \S_N$.
Since $h_N(\x)$ is a product of differences,
$\sum_{\sigma \in \S_N} {\rm sgn}(\sigma) \1_{\W_N^{\rm A}}
(\sigma(\B(t)))|h_N(\B(t))|
=\sum_{\sigma \in \S_N} \1_{\W_N^{\rm A}}(\sigma(\B(t)))
h_N(\B(t))=h_N(\B(t))$, and then (\ref{eqn:Asp2_1}) is 
simply written as
\begin{equation}
\E^{\xi}[g(\X(t))]
=\rE^{\x} \left[ g(\B(t))
\frac{h_N(\B(t))}{h_N(\x)} \right].
\label{eqn:Asp2_3}
\end{equation}

In the LHS of (\ref{eqn:Asp2_3}) we should note that
the Dyson model is an interacting particle system
such that between any pair of particles a long-range
repulsive force acts.
Strength of the two-body repulsive force is
exactly proportional to the inverse of distance
between two particles and thus it diverges
as the distance goes to zero.
By this strong repulsion, any collision of particles
is prevented.
On the other hand, in the RHS of (\ref{eqn:Asp2_3}),
independent BMs are considered.
When we calculate the expectation of a function $g$
of them, however, we have to put extra weight
$h_N(\B(t))/h_N(\x)$ to their paths.
Since if $|B_j(t)-B_i(t)| \to 0$ for any
$i \not= j$, this weight becomes zero, again any collision
of particle is prevented.
An important point of the Karlin-McGregor formula is that
this weight for paths is {\it signed},
{\it i.e.}, it can be positive and negative.
Therefore, all particle configurations 
realized by intersections of paths in the
1+1 spatio-temporal plane are completely cancelled.

\subsection{Complex BM representation}

Now we consider complexification of the expression 
(\ref{eqn:Asp2_3}) \cite{KT10b}.
For each $B_i^{x_i}(t), i \in \I$,
we introduce an independent one-dimensional BM
starting from the origin, $\widetilde{B}_i(t)$,
and define a complex BM as
$Z_i(t)=B_i(t)+\sqrt{-1} \widetilde{B}_i(t), i \in \I$.
If we write the expectation with respect to 
$\{\widetilde{B}_i(t)\}_{i \in \I}$
as $\widetilde{\rE}[\, \cdot \,]$
and define $\bE^{\x}=\rE^{\x} \otimes \widetilde{\rE}$,
we can confirm that the RHS of (\ref{eqn:Asp2_3}) can be
replaced by
\begin{equation}
\bE^{\x}  \left[g(\B(t))
\frac{h_N(\bZ(t))}{h_N(\x)} \right],
\label{eqn:cBM1}
\end{equation}
where $\bZ(t)=(Z_i(t))_{i \in \I}$.

A key lemma of our theory \cite{KT10b} is the following
identity;
for $\x=(x_1, \dots, x_N) \in \W_N^{\rm A}$,
$\z=(z_1, \dots, z_N) \in \C^{N}$,
$$
\frac{h_N(\z)}{h_N(\x)}
=\det_{1 \leq i, j \leq N}
\Big[ \Phi_{\xi}^{x_i}(z_j) \Big],
$$
where
\begin{equation}
\Phi_{\xi}^{u}(z)=\prod_{x \in \supp \xi \cap \{u\}^{\rm c}}
\left( 1- \frac{z-u}{x-u} \right)
\label{eqn:entire1}
\end{equation}
for $\xi=\sum_{i \in \I} \delta_{x_i} \in \mM_0$ and
$\supp \xi \equiv \{x \in \R: \xi(\{x\}) > 0\}$.
The function $\Phi_{\xi}^{u}(z)$ has an 
expression of the {\it Weierstrass canonical product with genus zero}.
Then, it is an {\it entire function} with zeros at
$\supp \xi \cap \{u\}^{\rm c}$ 
(see, for example, \cite{Lev96,Nog98}).

If we apply this identity to (\ref{eqn:cBM1}),
we have quantities $\Phi_{\xi}^{x_i}(Z_j(t)), i, j \in \I$,
which are conformal transforms of independent complex
BMs, $Z_j(t), j \in \I$.
Since complex BM is conformal invariant, each
$\Phi_{\xi}^{x_i}(Z_j(t))$ is a time change of
a complex BM, $Z_j(\cdot)$.
Then the average is conserved,
\begin{equation}
\bE^{x}[\Phi_{\xi}^{u}(Z_j(t))]
=\bE^{x}[\Phi_{\xi}^{u}(Z_j(T))],
\quad 0 \leq ^{\forall} t \leq T < \infty,
\label{eqn:martingale}
\end{equation}
that is, $\{\Phi_{\xi}^{u}(Z_j(t))\}_{j \in \I}$
are independent {\it conformal local martingales}
(see, for example, Section V.2 of \cite{RY98}).

Let $0 < t < T < \infty$. Then for any
$\F(t)$-measurable function $F$ in the continuous
path space $\c([0, \infty) \to \mM)$, we have
the equality
\begin{equation}
\E^{\xi} [F(\Xi(\cdot))]
=\bE^{\x} \left[ F\left( \sum_{i \in \I} \delta_{B_i(\cdot)} \right)
\det_{i, j \in \I} \Big[
\Phi_{\xi}^{x_i}(Z_j(T)) \Big] \right].
\label{eqn:cBM3}
\end{equation}
Now it is claimed that any observables of the Dyson model
is calculated by a system of independent complex BMs,
whose paths are weighted by a multivariate complex function
$\det_{i, j \in \I}[\Phi_{\xi}^{x_i}(Z_j(T))]$,
which is a conformal local martingale.
We call (\ref{eqn:cBM3}) the `complex BM representation'
of the Dyson model \cite{KT10b}.

\subsection{Determinantal process with an infinite
number of particles}

For a configuration $\xi \in \mM$, 
we write the restriction of configuration 
in $A \subset \R$ as
$(\xi\cap A) (\cdot)=\sum_{i \in \I : x_i \in  A}
\delta_{x_i}(\cdot)$, 
a shift of configuration by $u \in \R$ as
$\tau_u \xi(\cdot)=\sum_{i \in \I} \delta_{x_i+u}(\cdot)$,
and a square of configuration as 
$\xi^{\langle 2 \rangle}(\cdot)
=\sum_{i \in \I} \delta_{x_i^2}(\cdot)$, respectively.
Let $\c_0(\R)$ be the set of all continuous
real-valued functions with compact supports.

For any integer $M \in \N$,
a sequence of times
$\t=\{t_1,t_2,\dots,t_M \}$ with $0 < t_1 < \cdots < t_M < 
T < \infty$,
and a sequence of functions
$\f=(f_{t_1},f_{t_2},\dots,f_{t_M}) \in \c_0(\R)^M$,
the {\it moment generating function of multitime distribution}
of $(\Xi(t), \P^{\xi})$ is defined by
\begin{equation}
{\Psi}^{\xi}_{\t}[\f]
\equiv \E_{\xi} \left[\exp \left\{ \sum_{m=1}^{M} 
\int_{\R} f_{t_m}(x) \Xi(t_m, dx) \right\} \right].
\label{eqn:GF}
\end{equation}

From the fact that $\{\Phi_{\xi}^{u}(Z_j(t))\}_{j \in \I}$
are independent conformal local martingale
with the property (\ref{eqn:martingale}),
\begin{eqnarray}
\bE^{\x} \Big[ \Phi_{\xi}^{x_i}(Z_j(T)) \Big]
&=& \bE^{\x} \Big[ \Phi_{\xi}^{x_i}(Z_j(0)) \Big]
\nonumber\\
&=& \Phi_{\xi}^{x_i}(x_j) = \delta_{ij},
\quad ^{\forall}T \geq 0.
\nonumber
\end{eqnarray}
Moreover, we can see that
the following property holds for the 
determinantal weight in the complex BM representation
(\ref{eqn:cBM3}).
Let $\I'$ be a subset of the index set $\I$
and assume that a function $f$ depends on 
$B_i(t), i \in \I'$, but does not
on $B_j(t), j \in \I \setminus \I', 0 < t < T < \infty$.
Then
\begin{equation}
\bE^{\x} \left[ f(\{B_i(t)\}_{i \in \I'})
\det_{i, j \in \I} [ \Phi_{\xi}^{x_i}(Z_j(T))] \right]
= \bE^{\x'} \left[ f(\{B_i(t)\}_{i \in \I'})
\det_{i, j \in \I'} [ \Phi_{\xi}^{x_i}(Z_j(t))] \right],
\label{eqn:reduction}
\end{equation}
where $\bE^{\x'}[ \, \cdot \,]$ denotes
the expectation with respect to $\{B_i^{x_i}(t) \}_{i \in \I'}$.
Let $\xi \in \mM_0$ and
\begin{eqnarray}
\mbK^{\xi}(s, x; t, y) &=&
\int_{\R}\xi(dv)p_{s}(x|v)
\int_{\R} dw \ p_{t}(w|0) \Phi_{\xi}^{v}(y+\sqrt{-1}w)
\nonumber\\
&& -\1(s>t)p_{s-t}(x|y),
\quad (s,t) \in (0, \infty)^2, \quad
(x, y) \in \R^2,
\label{eqn:K1}
\end{eqnarray}
where $\1 (\omega)$ is 
the indicator function of a condition $\omega$;
$\1(\omega)=1$ if $\omega$ is satisfied 
and $\1(\omega)=0$ otherwise.
Then, using this `reducibility' (\ref{eqn:reduction}), 
we have proved that (\ref{eqn:GF}) 
is given by a {\it Fredholm determinant}
\begin{equation}
{\Psi}^{\xi}_{\t}[\f]
=\mathop{\Det}_
{\substack
{(s,t)\in \t^2, \\
(x,y)\in \R^2}
}
 \Big[\delta_{st} \delta(x-y)
+ \mbK^{\xi}(s,x;t,y) \chi_{t}(y) \Big],
\label{eqn:Fred}
\end{equation}
where
$\chi_{t_m}(\cdot)=e^{f_{t_m}(\cdot)}-1$, $1\leq m \leq M$.
We call $\mbK^{\xi}$ a {\it correlation kernel}.

For $L>0, \alpha>0$ and $\xi\in\mM$ we put
$$
M(\xi, L)=\int_{[-L,L]\setminus\{0\}} \frac{\xi(dx)}{x},
\qquad
M_\alpha(\xi, L)
=\left( \int_{[-L,L]\setminus\{0\}} 
\frac{\xi(dx)}{|x|^\alpha}\right)^{1/\alpha},
$$
and
$M(\xi) = \lim_{L\to\infty}M(\xi, L)$,
$M_\alpha(\xi)= \lim_{L\to\infty}M_\alpha(\xi, L)$,
if the limits finitely exist.
We have introduced the following conditions for
initial configurations $\xi \in \mM$ \cite{KT10}:

\vskip 3mm
\noindent ({\bf C.1})
there exists $C_0 > 0$ such that
$|M(\xi,L)|  < C_0$, $L>0$,

\vskip 3mm

\noindent ({\bf C.2}) (i) 
there exist $\alpha\in (1,2)$ and $C_1>0$ such that
$
M_\alpha(\xi) \le C_1,
$ \\
\noindent (ii) 
there exist $\beta >0$ and $C_2 >0$ such that
$$
M_1(\tau_{-a^2} \xi^{\langle 2 \rangle}) \leq C_2
(\max\{|a|, 1\})^{-\beta}
\quad \forall a \in \supp \xi.
$$
\vskip 3mm

\noindent
It was shown that, 
if $\xi\in\mM_0$ satisfies the 
conditions $({\bf C.1})$ and $({\bf C.2})$, then
for $a\in\R$ and $z\in\C$, 
$\displaystyle{\Phi_{\xi}^{a}(z)\equiv \lim_{L\to\infty}
\Phi_{\xi \cap [a-L, a+L]}^{a}(z)}$
finitely exists, and
$$
|\Phi_{\xi}^{a}(z)|\le C \exp \bigg\{c(|a|^\theta +|z|^\theta) \bigg\}
\left| \frac{z}{a}\right|^{\xi(\{0\})} \left|\frac{a}{a-z} \right|,
\quad a\in \supp \xi, \ z \in\C,
$$
for some $c, C>0$ and $\theta\in 
(\max\{\alpha, (2-\beta)\},2)$, 
which are determined by the constants $C_0, C_1, C_2$
and the indices $\alpha, \beta$ 
in the conditions \cite{KT10}.
Then even if $\xi(\R)=\infty$, under the conditions
$({\bf C.1})$ and $({\bf C.2})$, $\mbK^{\xi}$ given by 
(\ref{eqn:K1}) 
is well-defined as a correlation kernel and
dynamics of the Dyson model with an infinite number of
particles $(\Xi(t), \P^{\xi})$ exists \cite{KT10}.
We note that in the case that $\xi \in \mM_0$ satisfies the conditions 
$({\bf C.1})$ and $({\bf C.2})$
with constants $C_0, C_1, C_2$ and indices $\alpha$ and $\beta$,
then $\xi \cap [-L, L], \forall L > 0$ does as well.
Then we can obtain the convergence of moment generating functions
${\Psi}^{\xi \cap [-L, L]}_{\t}[\f] \to
{\Psi}^{\xi}_{\t}[\f]$ as $L \to \infty$, 
which implies the convergence of
the probability measures $\P^{\xi\cap [-L, L]}
\to \P^{\xi}$ in $L \to \infty$
in the sense of finite dimensional distributions \cite{KT10}.

By definition of Fredholm determinant, 
the moment generating function (\ref{eqn:Fred})
can be expanded with respect 
to $\chi_{t_m}(\cdot), 1 \leq m \leq M$, as
$$
{\Psi}^{\xi}_{\t}[\f]
=\sum_
{\substack
{N_m \geq 0, \\ 1 \leq m \leq M} }
\int_{\prod_{m=1}^{M} \W^{\rm A}_{N_{m}}}
\prod_{m=1}^{M} \left\{ d \x_{N_m}^{(m)}
\prod_{i=1}^{N_{m}} 
\chi_{t_m} \Big(x_{i}^{(m)} \Big) \right\}
\rho^{\xi} 
\Big( t_{1}, \x^{(1)}_{N_1}; \dots ; t_{M}, \x^{(M)}_{N_M} \Big),
$$
with
$$
\rho^{\xi} \Big(t_1,\x^{(1)}_{N_1}; \dots;t_M,\x^{(M)}_{N_M} \Big) 
=\det_{\substack
{1 \leq i \leq N_{m}, 1 \leq j \leq N_{n}, \\
1 \leq m, n \leq M}
}
\Bigg[
\mbK^{\xi}(t_m, x_{i}^{(m)}; t_n, x_{j}^{(n)} )
\Bigg],
$$
where $\x^{(m)}_{N_m}$ denotes
$(x^{(m)}_1, \dots, x^{(m)}_{N_m})$
and $d\x^{(m)}_{N_m}= \prod_{i=1}^{N_m} dx^{(m)}_i$,
$1 \leq m \leq M$.
The functions $\rho^{\xi}$'s are called {\it multitime 
correlation functions}, and ${\Psi}^{\xi}_{\t}[\f]$ 
can be regarded as 
a generating function of them.
In general, when the moment generating function for
the multitime distribution is given by
a Fredholm determinant,
all the spatio-temporal correlation functions
$\rho^{\xi}$ are given by determinants of matrices,
whose entries are special values of 
a continuous function $\mbK^{\xi}$,
and then 
the process is said to be 
{\it determinantal} \cite{NKT03,KNT04,KT07b}.
(Therefore, $\mbK^{\xi}$ is called a
correlation kernel.)
The results by Eynard and Mehta reported in \cite{EM98}
for a multi-layer matrix model can be regarded as the theorem that
the Dyson model is determinantal for the special
initial configuration $\xi =N \delta_0$, {\it i.e.},
all particles are put at the origin, 
for any $N \in \N$.
The correlation kernel is expressed in this case by using
the Hermite orthogonal polynomials \cite{NF98}.
The present author and H. Tanemura proved that, 
for any fixed initial configuration $\xi \in \mM$
with $\xi(\R) \in \N$,
the Dyson model $(\Xi(t), \P^{\xi})$
is determinantal,
in which the correlation kernel is given by
\begin{eqnarray}
&& \mbK^{\xi}(s,x; t,y)
= \frac{1}{2 \pi \sqrt{-1}} 
\oint_{\Gamma(\xi)} dz \, p_{s}(x|z)
\int_{\R} d w \, p_{t}(-\sqrt{-1} y|w)
\nonumber\\
&& \qquad \qquad \times
\frac{1}{\sqrt{-1} w-z} 
\prod_{x' \in \supp \xi}
\left( 1- \frac{\sqrt{-1}w-z}{x'-z} \right)
-{\bf 1}(s>t) p_{s-t}(x|y),
\label{eqn:K2}
\end{eqnarray}
where $\Gamma(\xi)$ is a closed contour on the
complex plane $\C$ encircling the points in 
$\supp \xi$ on the real line $\R$
once in the positive direction \cite{KT10}.
When $\xi \in \mM_0$, (\ref{eqn:K2}) becomes
(\ref{eqn:K1}) by performing the Cauchy integrals.

We note that
$\Xi_{\t} \equiv \sum_{t \in \t} \delta_t \otimes \Xi(t)$
is a {\it determinantal point process}
(or {\it Fermion point process}) on 
the spatio-temporal field $\t \times \R$
with an operator ${\cal K}$ given by
${\cal K}f(s,x)=\sum_{t \in \t} \int_{\R} dy \,
\mbK(s,x;t, y) f(t,y)$
for $f(t, \cdot) \in {\rm C}_0(\R), t \in \t$.
When ${\cal K}$ is symmetric, Soshnikov \cite{Sos00}
and Shirai and Takahashi \cite{ST03} gave sufficient
conditions for $\mbK^{\xi}$ to be a correlation kernel
of a determinantal point process 
(see also \cite{HKPV09}).
Such conditions are not known for asymmetric cases.
The present correlation kernels
(\ref{eqn:K1}) and (\ref{eqn:K2}) are asymmetric cases,
$\mbK^{\xi}(s,x;t,y) \not= \mbK^{\xi}(t,y;s,x)$
by the second terms 
$-\1(s>t) p_{s-t}(x|y)$.
Such form of asymmetric correlation kernels
is said to be of the Eynard-Mehta type \cite{BR05,KT10b}.

From the view point of statistical physics,
such asymmetry is useful to describe 
{\it nonequilibrium systems} developing in time.
In order to demonstrate it, we have studied
the following {\it relaxation phenomenon}
of the Dyson model with an infinite number of particles.

\begin{figure}
\begin{center}
\includegraphics[width=0.7\linewidth]{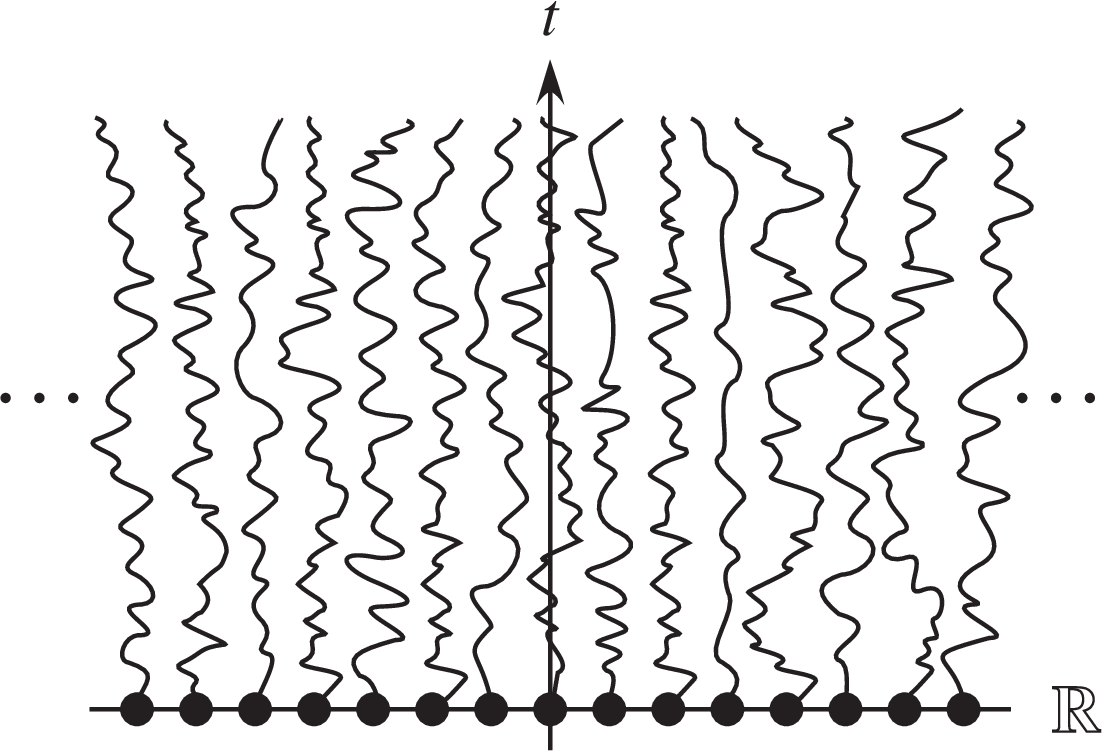}
\end{center}
\caption{\small
Consider the Dyson model starting from 
the configuration in which every point of the integers
$\Z$ is occupied by one particle.
This nonequilibrium determinantal process shows 
a relaxation phenomenon to the stationary state $\mu_{\sin}$.}
\label{fig:Dyson}
\end{figure}

We consider the configuration in which every point of the integers
$\Z$ is occupied by one particle,
$$
\xi^{\Z}(\cdot)
=\sum_{i \in \Z} \delta_{i}(\cdot).
$$
See Fig.\ref{fig:Dyson}.
It can be confirmed that $\xi^{\Z}$ satisfies the conditions
{\bf (C.1)} and {\bf (C.2)} and thus the Dyson model
starting from $\xi^{\Z}$, $(\Xi(t), \P^{\xi^{\Z}})$, 
is well-defined as a determinantal process
with an infinite number of particles.

As a matter of fact, we have shown that the correlation 
kernel is given by \cite{KT10}
\begin{eqnarray}
&& \mbK^{\xi^{\Z}}(s, x; t, y)
= {\bf K}_{\sin}(s, x; t, y) \nonumber\\
&& \qquad \qquad 
+ \frac{1}{2 \pi} \int_{|k| \leq \pi} dk \,
e^{k^2(t-s)/2+\sqrt{-1} k (y-x)}
\Big\{ \vartheta_3(x- \sqrt{-1} k s, 2 \pi \sqrt{-1} s) -1 \Big\}
\nonumber\\
&& \qquad = {\bf K}_{\sin}(s, x; t, y) \nonumber\\
&& \qquad \qquad + \sum_{n \in \Z \setminus \{0\}}
e^{2 \pi \sqrt{-1} x n - 2 \pi^2 s n^2}
\int_{0}^{1} du \,
e^{\pi^2 u^2 (t-s)/2}
\cos \Big[ \pi u 
\{ (y-x) - 2 \pi \sqrt{-1} s n \} \Big],
\nonumber
\end{eqnarray}
$(s, t) \in [0, \infty)^2, (x,y) \in \R^2$,
where 
\begin{eqnarray}
{\bf K}_{\sin}(s, x; t, y) &=&
\frac{1}{2 \pi} \int_{|k| \leq \pi} dk \,
e^{k^2(t-s)/2 + \sqrt{-1} k (y-x)}
- {\bf 1}(s>t) p_{s-t}(x|y) \nonumber\\
&=& \left\{ \begin{array}{ll} 
\displaystyle{
\int_{0}^{1} du \, e^{\pi^2 u^2 (t-s)/2} 
\cos \{ \pi u (y-x)\} }
& \mbox{if $t>s $} \cr
& \cr
K_{\sin}(x,y)
& \mbox{if $t=s$} \cr
& \cr
\displaystyle{
- \int_{1}^{\infty} du \, 
e^{\pi^2 u^2 (t-s)/2} \cos \{ \pi u (y-x) \} }
& \mbox{if $t<s$},
\end{array} \right.
\nonumber
\end{eqnarray}
with
\begin{equation}
K_{\sin}(x,y)= 
\frac{1}{2 \pi} \int_{|k| \leq \pi} dk \,
e^{\sqrt{-1} k(y-x)}
= \frac{\sin \{ \pi(y-x) \} }{\pi (y-x)},
\quad x, y \in \R,
\label{eqn:sine_kernel0}
\end{equation}
and $\vartheta_3$ is a version of the {\it Jacobi theta function}
defined by
\begin{equation}
\vartheta_3(v, \tau) 
= \sum_{n \in \Z} e^{2 \pi \sqrt{-1} v n+\pi \sqrt{-1} \tau n^2},
\quad \Im \tau > 0.
\label{eqn:theta3}
\end{equation}
By this explicit expression, we can see that
\begin{equation}
\lim_{u \to \infty}
\mbK^{\xi^{\Z}}(u+s,x; u+t, y)={\bf K}_{\sin}(s,x;t,y).
\label{eqn:conv}
\end{equation}
The correlation kernel ${\bf K}_{\sin}(s,x;t,y)$
surviving in the long-term limit
is called the {\it extended sine kernel}.
It is symmetric,
${\bf K}_{\sin}(s,x;t,y)={\bf K}_{\sin}(t,y;s,x)$,
and the determinantal process with the
correlation kernel ${\bf K}_{\sin}(s,x;t,y)$ is 
an {\it equilibrium dynamics}.
It is time-reversal
with respect to the determinantal point process
$\mu_{\sin}$, in which any spatial correlation
function is given by a determinant with
the correlation kernel (\ref{eqn:sine_kernel0})
called the {\it sine kernel}.
This stationary state $\mu_{\sin}$ is 
a scaling limit (called the bulk-scaling-limit)
of the eigenvalue distribution of random matrices
in GUE \cite{Meh04,For10}. 
See \cite{Spo87,Osa96,NF98,KT07b}.

The theory of entire functions discusses
the relations between the growth of an entire function
and the distribution of its zeros \cite{Lev96,Nog98}.
Here we set distributions of zeros as
initial configurations satisfying some conditions
and control the behavior of particles at infinity
to realize nonequilibrium dynamics of 
long-rang interacting infinite-particle systems.
Systematic study of determinantal processes with
infinite number of particles exhibiting relaxation phenomena
is now in progress \cite{KT09,KT10,KT11,KT10b}
\footnote{
It will be interesting to discuss intrinsic relations
between the above mentioned relaxation phenomenon
to $\mu_{\sin}$ and the 
{\it interpolation/sampling theorem} 
of Whittaker and others \cite{BFHSSS10},
which represents a function $f$ 
in the cardinal sampling series
$$
f(z)=\frac{\sin \pi z}{\pi}
\sum_{n \in \Z} f(n) \frac{(-1)^n}{z-n}.
$$
}.

\SSC{Related Topics}

At the end of this manuscript,
I briefly introduce related topics,
which we are interested in.

\subsection{Extreme value distributions of
noncolliding diffusion processes}

As explained in Section 1.2, BES$^{(3)}$ is the 
conditional BM to stay positive.
When we impose an additional condition such that
it starts from the origin at time $t=0$
and return to the origin at time $t=1$,
the process is called the
{\it three-dimensional Bessel bridge}
with duration 1,
which is here denoted by $Y(t), 0 \leq t \leq 1$.

\begin{figure}
\begin{center}
\includegraphics[width=0.5\linewidth]{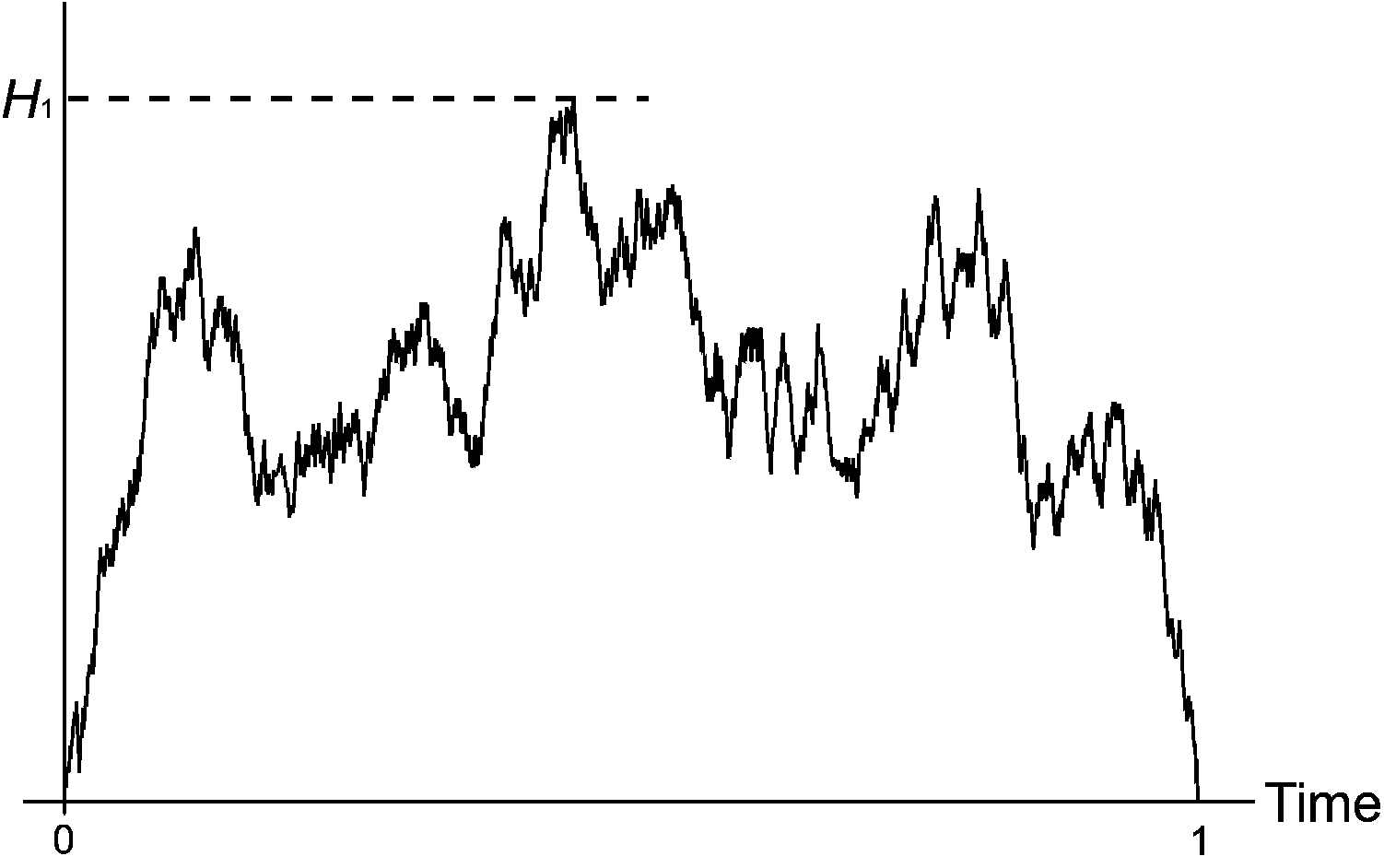}
\end{center}
\caption{\small
Let $H_1$ be the maximum of the three-dimensional
Bessel bridge with duration 1.
The expectation $\rE[H_1^s]$ is expressed by
using the value of Riemann zeta function at 
$s \in \C$.}
\label{fig:extreme}
\end{figure}

Here we consider the maximum of $Y(t)$,
$$
H_1 = \max_{0 < t < 1} Y(t).
$$
See Fig.\ref{fig:extreme}
We can show that the distribution of $H_1$ is 
described as \cite{KIK08b}
\begin{equation}
\rP(H_1 \leq h)
=-\frac{1}{2} \sum_{n \in \Z} H_2(\sqrt{2} n h)
e^{-2 h^2 n^2},
\label{eqn:max2}
\end{equation}
where $H_{i}(x)$ is the $i$-th Hermite polynomial
\begin{equation}
H_{i}(x)=i ! \sum_{j=0}^{[i/2]}
\frac{(-1)^{j} (2x)^{i-2j}}
{j ! (i - 2j)!},
\quad i \in \{0,1,2, \cdots\}
\label{eqn:Hermite1}
\end{equation}
with $[a]$=the greatest integer that is not 
greater than $a \in \R$.
By the equation 
$\rP(H_1 \leq h) = \int_0^h p_1(u) du$,
the probability density function $p_1(u)$ of $H_1$
is defined, and $s$-th moment of 
the random variable $H_1$ is calculated by
$$
\rE[ H_1^s ] = \int_0^{\infty} dh \,
h^s p_1(s).
$$
As discussed by Biane, Pitman, and Yor \cite{BPY01},
if we set $\rE[H_1^s]=2(\pi/2)^{s/2} \xi(s)$,
the following equality is established,
\begin{eqnarray}
\xi(s) &=& \frac{1}{2} s(s-1) \pi^{-s/2} 
\Gamma(s/2) \zeta(s) \nonumber\\
&=& \frac{1}{2}+\frac{1}{4} s(s-1)
\int_1^{\infty} du \, (u^{s/2-1}+u^{(1-s)/2-1})
(\vartheta_3(0, \sqrt{-1} u)-1),
\label{eqn:xis}
\end{eqnarray}
where $\Gamma(z)$ is the gamma function,
$\zeta(z)$ is the {\it Riemann zeta function},
\begin{equation}
\zeta(z)=\sum_{n=1}^{\infty} \frac{1}{n^z},
\quad \Re z > 1,
\label{eqn:zeta}
\end{equation}
and $\vartheta_3(v, \tau)$ is given by
(\ref{eqn:theta3}).
See also Chapter 11 in \cite{Yor97}.

Then we consider an $N$-tuples of three-dimensional
Bessel bridges 
conditioned never to collide with each other,
$\Y^{(N)}(t)=(Y^{(N)}_1(t), \dots, Y^{(N)}_N(t))$,
$$
0< Y^{(N)}_1 < Y^{(N)}_2 < \cdots
< Y^{(N)}_N,
\quad 0 < t < 1,
$$
with $\Y^{(N)}(0)=\Y^{(N)}(1)=\0$.
It is found that (\ref{eqn:max2}) is generalized for
$$
H_N=\max_{0 < t < 1} Y^{(N)}_N
$$
as \cite{SMCRF08,KIK08b}
\begin{equation}
\rP(H_N \leq h)
=\frac{(-1)^{N}}{2^{N^2} \prod_{i=1}^{N} \Gamma(2i)}
\det_{1 \leq i, j \leq N}
\left[ \sum_{n \in \Z}
H_{2(i+j-1)}(\sqrt{2} n h) e^{-2 n^2 h^2} \right].
\label{eqn:max5}
\end{equation}
Extensive study of extreme value distributions of
noncolliding processes has been reported
in \cite{Ful07,KIK08a,Fei07,SMCRF08,Fei09,BFPSW09,NM09,RS10,IK11}.
See also \cite{KO01,TW07,KT07a,KMW09}.

Recently, Forrester, Majumdar, and Schehr
clarified an equality between the distribution function
$\rP(H_N \leq h)$ and the partition function
of the 2$d$ {\it Yang-Mills theory} on a sphere with
a gauge group ${\rm Sp}(2N)$ \cite{FMS11}.
Moreover, by using this equivalence,
they proved that in a scaling limit of $\rP(H_N \leq h)$,
the {\it Tracy-Widom distribution} of GOE
of random matrices \cite{TW94,TW96} is derived.

\subsection{Characteristic polynomials of random matrices}

Let $\cH$ be an $N \times N$ Hermitian random matrix
in the GUE. The characteristic polynomial of
a variable $\alpha$ is then given by
$$
P(\alpha)=\det ( \alpha I - \cH),
$$
where $I$ is the $N \times N$ unit matrix.
In the connection with the Riemann zeta function
(\ref{eqn:zeta}),
statistical property of $P(\alpha)$ has been
studied \cite{KS00a,KS00b,HKO00,BH00}.
Here we consider the ensemble average of
$m$-product of characteristic polynomials
\begin{equation}
M_{\rm GUE}(m, \valpha; N, \sigma^2)
= \left\langle
\prod_{n=1}^{m} P(\alpha_n) 
\right\rangle_{{\rm GUE}(N, \sigma^2)}
= \left\langle
\prod_{n=1}^{m} \prod_{i=1}^{N} (\alpha_n -\lambda_i) 
\right\rangle_{{\rm GUE}(N, \sigma^2)}, 
\label{eqn:chara2}
\end{equation}
where $\left\langle \, \cdot \,
\right\rangle_{{\rm GUE}(N, \sigma^2)}$
denotes the ensemble average in the GUE
with variance $\sigma^2$ of $N \times N$
Hermitian matrices $\{ \cH \}$,
whose eigenvalues are $\{(\lambda_i)_{i=1}^N\}$.
The probability density function of eigenvalues of
${\cH}$ in GUE with variance $\sigma^2$ is given by
\begin{equation}
\mu_{N, \sigma^2}(\xi)
=\frac{\sigma^{-N^2}}{C_N}
\exp \left( - \frac{|\x|^2}{2 \sigma^2} \right)
h_N(\x)^2,
\label{eqn:GUE1}
\end{equation}
$\xi=\sum_{i=1}^{N} \delta_{x_i} \in \mM, 
x_1 \leq x_2 \leq \cdots \leq x_N$,
where $C_N=(2 \pi)^{N/2} \prod_{i=1}^{N} \Gamma(i)$, 
$|\x|^2=\sum_{i=1}^{N} x_i^2$,
and $h_N(\x)$ is given by (\ref{eqn:hN2}).
Then if the expectation of a measurable function $F$ of
a random variable $\Xi \in \mM$ with respect to
(\ref{eqn:GUE1}) is written as
$$
\bE_{N, \sigma^2} [F(\Xi)]
=\int_{\W_N^{\rm A}} F(\xi) \mu_{N, \sigma^2}(\xi) d \x
$$
with setting $\xi=\sum_{i=1}^{N} \delta_{x_i},
\x=(x_1, \dots, x_N)$,
where $d \x=\prod_{i=1}^N dx_i$,
(\ref{eqn:chara2}) is written as
$$
M_{\rm GUE}(m, \valpha; N, \sigma^2)
= \bE_{N, \sigma^2} \left[
\prod_{n=1}^{m} \prod_{X \in \Xi}(\alpha_n-X) \right].
$$

In Section 3.3, we showed that the Dyson model $(\Xi(t), \P^{\xi})$
is a determinantal process with the correlation kernel
(\ref{eqn:K1}) for any fixed initial configuration
$\xi \in \mM_0$ if $\xi(\R)=N \in \N$.
Here we consider the situation such that the initial
configuration $\xi$ is distributed according to
(\ref{eqn:GUE1}).
Note that by the term $h_N(\x)^2$ in (\ref{eqn:GUE1}),
the GUE eigenvalue distribution is in $\mM_0$
w.p.1.

By using Aspect 2 of the Dyson model,
we can prove that this process, denoted by
$(\Xi(t), \P^{\mu_{N, \sigma^2}})$, is
equivalent with the time shift $t \to t+\sigma^2$
of the Dyson model starting from the configuration $N \delta_0$
({\it i.e.}, all $N$ particles are put at the origin).
That is, the equality
\begin{equation}
(\Xi(t), \P^{\mu_{N, \sigma^2}}) 
= (\Xi(t+\sigma^2), \P^{N \delta_0})
\label{eqn:chara4}
\end{equation}
holds for arbitrary $\sigma^2 >0$ in the sense of
finite dimensional distribution \cite{Kat11}.

This equivalence is highly nontrivial, since
even from its special consequence,
the following 
determinantal expression is derived for 
the ensemble average of product of 
characteristic polynomials; 
For any $N, n \in \N$, 
$\valpha=(\alpha_1, \alpha_2, \cdots, \alpha_{2n}) \in \C^{2n}$,
$\sigma^2>0$, 
\begin{eqnarray}
&&M_{\rm GUE}(2n, \valpha; N, \sigma^2)
=  
\frac{\gamma_{N,2n} \sigma^{n(2N+n)}}
{h_n(\alpha_1, \cdots, \alpha_n) h_n(\alpha_{n+1}, 
\cdots, \alpha_{2n})}
\nonumber\\
&& \quad \times
\det_{1 \leq i, j \leq n}
\left[ \frac{1}{\alpha_i-\alpha_{n+j}}
\left| \begin{array}{cc}
H_{N+n} (\alpha_i/\sqrt{2 \sigma^2} ) 
& H_{N+n} ( \alpha_{n+j}/\sqrt{2 \sigma^2} ) \cr
H_{N+n-1} ( \alpha_i/\sqrt{2 \sigma^2} ) 
& H_{N+n-1}( \alpha_{n+j}/\sqrt{2 \sigma^2} )
\end{array} \right| \right], 
\label{eqn:charaA1}
\end{eqnarray}
where
$$
\gamma_{N,2n}=
2^{-n(2N+2n-1)/2} \prod_{i=2}^{n}
\frac{(N+n-i)!}{(N+n-1)!},
$$
and $H_{i}(x)$ is the $i$-th 
Hermite polynomial given by (\ref{eqn:Hermite1}).

Moreover, in order to simplify the above expression, we can use
the following identity, which was recently
given by Ishikawa {\it et al.} \cite{IOTZ06} 
as a generalization of the Cauchy determinant
$$
\det_{1 \leq i, j \leq n}
\left( \frac{1}{x_i+y_j} \right)
= \frac{h_n(\x) h_n(\y)}
{\prod_{i=1}^{n} \prod_{j=1}^{n} (x_i+y_j)}.
$$
For $n \geq 2, \x=(x_1, \dots, x_n), 
\y=(y_1, \dots, y_n),
\a=(a_1, \dots, a_n), \b=(b_1, \dots, b_n) \in \C^n$,
\begin{eqnarray}
&& \det_{1 \leq i, j \leq n}
\left[ \frac{1}{y_j-x_i} 
\left| \begin{array}{cc}
1 & a_i \cr 1 & b_j 
\end{array} \right| \right]
\nonumber\\
&=& \frac{(-1)^{n(n-1)/2}}
{\prod_{i=1}^{n} \prod_{j=1}^{n}(y_j-x_i)}
\left| \begin{array}{cccccccc}
1 & x_1 & \cdots & x_1^{n-1} & a_1 & a_1 x_1 & \cdots a_1 x_1^{n-1} \cr
1 & x_2 & \cdots & x_2^{n-1} & a_2 & a_2 x_2 & \cdots a_2 x_2^{n-1} \cr
  &     & \cdots &           &     &         & \cdots               \cr
1 & x_n & \cdots & x_n^{n-1} & a_n & a_n x_n & \cdots a_n x_n^{n-1} \cr
1 & y_1 & \cdots & y_1^{n-1} & b_1 & b_1 y_1 & \cdots b_1 y_1^{n-1} \cr
1 & y_2 & \cdots & y_2^{n-1} & b_2 & b_2 y_2 & \cdots b_2 y_2^{n-1} \cr
  &     & \cdots &           &     &         & \cdots               \cr
1 & y_n & \cdots & y_n^{n-1} & b_n & b_n y_n & \cdots b_n y_n^{n-1} \cr
\end{array} \right|.
\label{eqn:Ishikawa1}
\end{eqnarray}
Then we have the expression
\begin{equation}
M_{\rm GUE}(2n, \valpha; N, \sigma^2)= 
\frac{1}{h_{2n}(\valpha)}
\det_{1 \leq i, j \leq 2n}
\left[ \widehat{H}_{N+i-1}(\alpha_j; \sigma^2) \right],
\label{eqn:charaA2}
\end{equation}
where
$$
\widehat{H}_{i}(\alpha; \sigma^2)
\equiv \left(\frac{\sigma^2}{2} \right)^{i/2}
H_{i} \left(\frac{\alpha}{\sqrt{2 \sigma^2}}\right).
$$
This determinantal expression (\ref{eqn:charaA2})
is also obtained from the general formula given by
Br\'ezin and Hikami as Eq.(14) in \cite{BH00}.
See \cite{Del10} for recent development of
this topic.

\subsection{Fomin's determinant for loop-erased
random walks and its scaling limit}

We consider a network $\Gamma=(V,E,W)$,
where $V=\{v_i\}$ and $E=\{e_i\}$ are
sets of vertices and of edges
of an undirected planar lattice, respectively,
and $W=\{w(e)\}_{e \in E}$ is
a set of the weight functions of edges.
For $a, b \in V$, 
let $\pi$ be a walk given by
$$
\pi \, : \,
a=v_0 \,
{\stackrel{e_1}{\rightarrow}} \, 
v_1 \, 
{\stackrel{e_2}{\rightarrow}} \, 
v_2 \, 
{\stackrel{e_3}{\rightarrow}} \, 
\cdots
{\stackrel{e_m}{\rightarrow}} \, 
v_m =b \,
$$
where the length of walk is 
$|\pi|=m \in \N$ and, 
for each $0 \leq i \leq m-1$, 
$v_i$ and $v_{i+1}$ are nearest-neighboring vertices
in $V$ and $e_i \in E$ is the edge connecting
these two vertices. 
The weight of $\pi$ is given by 
$w(\pi)=\prod_{i=1}^{m} w(e_i)$.
For any two vertices of $a, b \in V$,
the Green's function of walks 
$\{\pi : a \rightarrow b\}$
is defined by
$$
W(a,b)=\sum_m \sum_{\pi: a \rightarrow b , |\pi|=m }
w(\pi).
$$
The matrix $W=(W(a,b))_{a, b \in V}$ is called
the {\it walk matrix} of the network $\Gamma$.

The loop-erased part of $\pi$,
denoted by $\LE(\pi)$, is defined recursively as follows.
If $\pi$ does not have self-intersections,
that is, all vertices $v_i, 0 \leq i \leq m$ are distinct, then
$\LE(\pi)=\pi$.
Otherwise, set $\LE(\pi)=\LE(\pi')$,
where $\pi'$ is obtained by removing the first loop
it makes.
The loop-erasing operator $\LE$ maps
arbitrary walks to {\it self-avoiding walks} (SAWs).
Note that the map is many-to-one.
For each SAW, $\zeta$, the weight $\widetilde{w}(\zeta)$
is given by
\begin{equation}
\widetilde{w}(\zeta)=
\sum_{\pi : \LE(\pi)=\zeta}
w(\pi).
\label{eqn:LEweight}
\end{equation}
We consider the statistical ensemble of
SAWs with the weight (\ref{eqn:LEweight})
and call it {\it loop-erased random walks} (LERWs) \cite{LL10}.

\begin{figure}
\begin{center}
\includegraphics[width=0.5\linewidth]{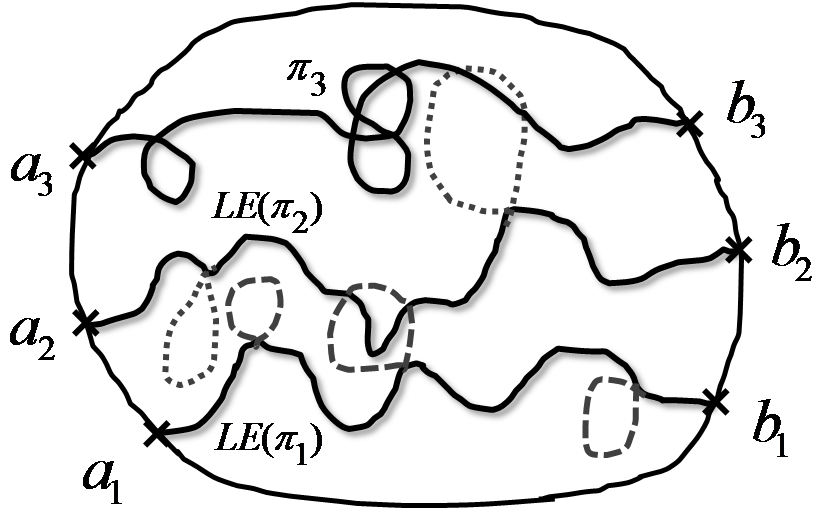}
\end{center}
\caption{\small
The situation $\LE(\pi_j) \cap \pi_3=\emptyset,
j=1,2$ is illustrated in a planar domain $D$,
where $A=(a_1, a_2, a_3)$ and $B=(b_1, b_2, b_3)$
are all boundary points of $\partial D$.
In this figure, $\LE(\pi_1)$ and $\LE(\pi_2)$ 
denoted by solid curves are
the loop-erased parts of the walks
$\pi_1: a_1 \rightarrow b_1$ and
$\pi_2: a_2 \rightarrow b_2$, respectively.
The third walk $\pi_3: a_3 \rightarrow b_3$
can be self-intersecting, but it does not
intersect with $\LE(\pi_1)$ nor $\LE(\pi_2)$.
}
\label{fig:LERW}
\end{figure}

Assume that $A=\{a_1, a_2, \dots, a_N\} \subset V$ and
$B=\{b_1, b_2, \dots, b_N \} \subset V$ are chosen
so that any walk from $a_i$ to $b_j$ intersects
any walk from $a_{i'}, i' > i$,
to $b_{j'}, j' < j$.
The weight of $N$-tuples of independent walks
$ a_1 \, {\stackrel{\pi_1}{\rightarrow}} \, b_1, 
\dots, 
a_N \, {\stackrel{\pi_N}{\rightarrow}} \, b_N $
is given by the product of $N$ weights
$\prod_{i=1}^N w(\pi_{i})$.
Then we consider $N$-tuples of walks
$(\pi_1, \pi_2, \dots, \pi_N)$ 
conditioned so that, for any
$1 \leq i < j \leq N$, the walk $\pi_j$
has no common vertices with the loop-erased part
of $\pi_i$;
\begin{equation}
\LE(\pi_i) \cap \pi_j=\emptyset, \quad
1 \leq i < j \leq N.
\label{eqn:LEW0}
\end{equation}
See Fig.\ref{fig:LERW}.
By definition, $\LE(\pi_j)$ is a part of $\pi_j$, and thus
nonintersection of any pair of loop-erased parts
is concluded from (\ref{eqn:LEW0});
$$
\LE(\pi_i) \cap \LE(\pi_j) = \emptyset,
\quad 1 \leq i < j \leq N.
$$

Fomin proved that total weight of $N$-tuples of walks
satisfying such a version of nonintersection condition
is given by the minor of
walk matrix,
$\det(W_{A,B}) \equiv
\det_{a \in A, b \in B}(W(a,b))$ \cite{For01}.
This minor is called {\it Fomin's determinant}
and Fomin's formula is expressed by
the equality \cite{For01,LL10}
\begin{equation}
\det(W_{A,B}) = 
\sum_{\LE(\pi_i) \cap \pi_j=\emptyset, \,
i<j} 
\prod_{k=1}^{N} w(\pi_{k}).
\label{eqn:Fomin}
\end{equation}

Kozdron and Lawler \cite{KL05}
consider continuum limit (the diffusion scaling limit) of 
Fomin's determinantal system of 
loop-erased random walks
in the complex plane $\C$,
where the initial and the final points
$A=\{a_i\}$ and $B=\{b_i\}$ of paths 
can be put on the boundaries of the domains
$\partial D$. 
By the diffusion scaling limit
each random walk will converge to 
a path of complex BM. 
We should note that, however, the characteristics of BM
look more similar to those of a surface than
those of a curve.
It implies that the Brownian path has loops on every scale
and then the loop-erasing procedure mentioned above
does not make sense for BM in the plane,
since we can not decide which loop is the first one.
Kozdron and Lawler proved explicitly, however, that
the continuum limit of Fomin's determinant of the
Green's functions of random walks
converges to that of the Green's functions of 
Brownian motions \cite{KL05}.
This will enable us to discuss 
{\it nonintersecting systems of 
loop-erased Brownian paths} in the sense of Fomin (\ref{eqn:LEW0}).
Moreover, Kozdron \cite{Koz09} showed that 
$2 \times 2$ Fomin's determinant
representing the event $\LE(\beta_1) \cap \beta_2=\emptyset$
for two complex Brownian paths $(\beta_1, \beta_2)$
is proportional to the probability that 
$\gamma_{{\rm SLE}^{(3)}} \cap \beta=\emptyset$,
where $\gamma_{{\rm SLE}^{(3)}}$ and $\beta$ denote 
the SLE$^{(3)}$ path
and a complex Brownian path, respectively.
On the other hand, Lawler and Werner gave a correct way
to add `Brownian loops' to an SLE$^{(3)}$ path
to obtain a complex Brownian path \cite{LW04}.
These results imply that the scaling limit of
loop-erased part of complex Brownian path
is described by the SLE$^{(3)}$ path,
as announced in Section 1.4.

Setting a sequence of chambers
in a planar domain, Sato and the present author observe
the first passage points at which $N$-tuples 
of complex Brownian paths
$(\beta_1, \dots, \beta_N)$ first enter each chamber, 
under the condition that the loop-erased parts
$(\LE(\beta_1), \dots, \LE(\beta_N))$ 
make a nonintersecting system in the domain
in the sense of Fomin (\ref{eqn:LEW0}) \cite{SK11}.
It is proved that the system of first passage points is
a determinantal point process in the planar domain,
in which the correlation kernel is of Eynard-Mehta type
\cite{SK11}.
Interpretation of this result
in terms of `mutually avoiding SLE paths' 
\cite{KL07,Dub06}
will be an interesting future problem.

\vskip 0.5cm
\begin{small}
\noindent{\bf Acknowledgements} \quad
The present author expresses his gratitude for 
Jun-ichi Matsuzawa for giving him such an opportunity to give
a talk at the Oka symposium.
This manuscript is based on the joint work with 
Hideki Tanemura, Taro Nagao, Naoki Kobayashi, 
Naoaki Komatsuda, Minami Izumi, and Makiko Sato.


\end{small}
\end{document}